\documentclass[a4paper]{amsart}
\usepackage{amsfonts}
\usepackage{amssymb}
\usepackage{euscript}
\usepackage{xspace}
\usepackage{enumerate}

\newtheorem{theorem}{Theorem}[section]
\newtheorem{proposition}[theorem]{Proposition}
\newtheorem{corollary}[theorem]{Corollary}
\newtheorem{lemma}[theorem]{Lemma}
\theoremstyle{definition}
\newtheorem{definition}[theorem]{Definition}

\newcommand{\blank}{{\llcorner\!\lrcorner}}
\newcommand{\Mp}[1]{${#1}$\hbox{-}\penalty60\hspace*{0pt}}
\newcommand{\Mpn}[1]{${#1}$\hbox{-}\penalty10000\hspace*{0pt}}

\newcommand{\Adj}{\mathbb{B}}        
\newcommand{\AdjG}{\mathbb{B}^G}     
\newcommand{\Comp}{\mathbb{K}}       
\newcommand{\Mat}{\mathbb{M}}        
\newcommand{\Mult}{\mathcal{M}}      

\newcommand{\KET}[1]{{\lvert#1\mathclose{\rangle\!\rangle}}}
\newcommand{\BRA}[1]{{\mathopen{\langle\!\langle}#1\rvert}}
\newcommand{\BRAKET}[2]{{\mathopen{\langle\!\langle}#1\mid#2\mathclose
{\rangle\!\rangle}}}

\newcommand{\ket}[1]{{\lvert#1\rangle}}
\newcommand{\bra}[1]{{\langle#1\rvert}}
\newcommand{\braket}[2]{{\langle#1\mid#2\rangle}}

\newcommand{\Fix}{\operatorname{Fix}} 
\newcommand{\Aut}{\operatorname{Aut}} 
\newcommand{\End}{\operatorname{End}} 

\newcommand{\si}{{\mathrm{si}}}      
\newcommand{\alg}{{\mathrm{alg}}}    
\newcommand{\cs}{{\mathrm{c}}}       

\newcommand{\C}{{\mathbb{C}}}        
\newcommand{\R}{{\mathbb{R}}}        
\newcommand{\Z}{{\mathbb{Z}}}        
\newcommand{\N}{{\mathbb{N}}}        
\newcommand{\exN}{\overline{\N}}     

\newcommand{\norm}[1]{\lVert#1\rVert} 

\newcommand{\ID}{\mathrm{id}}        
\newcommand{\defeq}{\mathrel{:=}}

\newcommand{\Cstar}{\Mpn{C^\ast}}
\newcommand{\Cred}{C^\ast_{\mathrm{r}}}

\newcommand{\Hils}{\EuScript{H}}     
\newcommand{\E}{\EuScript{E}}        
\newcommand{\F}{\EuScript{F}}        
\newcommand{\EL}{\EuScript{L}}       
\newcommand{\RC}{\EuScript{R}}       
\newcommand{\prid}{{\mathfrak{p}}}   

\begin{document}
\title[Generalized fixed point algebras]{Generalized fixed point algebras and
square-integrable group actions}
\author{Ralf Meyer}
\address{Mathematisches Institut\\
         Westf\"alische Wilhelms-Universit\"at M\"unster\\
         Einsteinstr.\ 62\\
         48149 M\"unster\\
         Germany
}
\email{rameyer@math.uni-muenster.de}
\urladdr{http://www.math.uni-muenster.de/u/rameyer}
\begin{abstract}
  We analzye Rieffel's construction of generalized fixed point algebras in the
  setting of group actions on Hilbert modules.  Let~$G$ be a locally compact
  group acting on a $C^\ast$-algebra~$B$.  We construct a Hilbert module~$F$
  over the reduced crossed product of $G$ and~$B$, using a pair $(E, R)$,
  where~$E$ is an equivariant Hilbert module over~$B$ and~$R$ is a dense
  subspace of~$E$ with certain properties.  The generalized fixed point
  algebra is the $C^\ast$-algebra of compact operators on~$F$.  Any Hilbert
  module over the reduced crossed product arises by this construction for a
  pair $(E, R)$ that is unique up to isomorphism.

  A necessary condition for the existence of~$R$ is that~$E$ be
  square-integrable.  The consideration of square-integrable representations
  of Abelian groups on Hilbert space shows that this condition is not
  sufficient and that different choices for~$R$ may yield different
  generalized fixed point algebras.

  If~$B$ is proper in Kasparov's sense, there is a unique~$R$ with the
  required properties.  Thus the generalized fixed point algebra only depends
  on~$E$.
\end{abstract}
\maketitle

\section{Introduction}
\label{sec:intro}

Let $(G, X, \alpha)$ be a dynamical system, consisting of a locally compact
group~$G$, a locally compact space~$X$, and a continuous left action $\alpha
\colon G \times X \to X$.  The action is called \emph{proper} iff for all
compact subsets $K, L \subseteq X$, the set of $g \in G$ with $\alpha_g (K)
\cap L \neq \emptyset$ is (relatively) compact.  Proper actions have many nice
properties.  For instance, the orbit space $G \backslash X$ is again a locally
compact space.  Rieffel~\cite{rieffel:proper} has initiated a program to
extend the notions of proper action and orbit space to noncommutative
dynamical systems, that is, group actions on \Cstar{}algebras.

Suppose that the group~$G$ is compact.  Then all actions of~$G$ on
\Cstar{}algebras are proper.  The role of the orbit space is played by the
\emph{fixed point algebra}
$$
A^G \defeq
\{ a \in A \mid \text{$\alpha_g (a) = a$ for all $g \in G$} \}.
$$
This is reasonable because $A^G \cong C_0 (G \backslash X)$ if $A = C_0 (X)$.

If~$G$ fails to be compact, there are several ways to define ``proper''
actions on \Cstar{}algebras.  The weakest reasonable notion is
\emph{square-integrability}.  It has interesting applications in equivariant
Kasparov theory~\cite{meyer:KKG} but is not enough to construct an ``orbit
space'', that is, a generalized fixed point algebra.  A slightly more
restrictive assumption is \emph{continuous square-integrability}, which is
exactly what is needed to construct a generalized fixed point algebra.
Another much more restrictive notion of properness is due to Kasparov (see
below).  To avoid a conflict of notation, we only use the word ``proper'' in
Kasparov's sense.

We are going to explain square-integrability and a variant of Rieffel's
construction of generalized fixed point algebras.  For both purposes, it is
very illuminating to allow group actions on Hilbert modules, not just on
\Cstar{}algebras.  Hilbert modules are more flexible because there are always
plenty of adjointable operators between them.  Let~$B$ be a
\Mpn{G}\Cstar{}algebra, let~$\E$ be a \Mpn{G}equivariant Hilbert module
over~$B$, and let $\xi \in \E$.  Denote the action of~$G$ on~$\E$ by~$\gamma$.
Define
\begin{xalignat}{2}
  \label{eq:def_BRA}
  \BRA {\xi} & \colon \E \to  C_b (G, B), &
  (\BRA {\xi} \eta) (g) &\defeq
  \braket {\gamma_g (\xi)} {\eta},
  \\
  \label{eq:def_KET}
  \KET {\xi} & \colon C_c (G, B) \to \E, &
  \KET {\xi} f & \defeq
  \int_G \gamma_g (\xi) \cdot f (g) \,dg.
\end{xalignat}
The operators $\KET {\xi}$ and $\BRA {\xi}$ are \Mpn{G}equivariant and adjoint
to each other with respect to the pairing $\braket {f_1} {f_2} \defeq \int_G
f_1 (g)^\ast f_2 (g) \,dg$ between $C_c (G, B)$ and $C_b (G, B)$.

We call~$\xi$ \emph{square-integrable} iff $\BRA {\xi} \eta \in L^2 (G, B)$
for all $\eta \in \E$.  Let $\E_\si \subseteq \E$ be the subspace of
square-integrable elements.  If $\xi \in \E_\si$, then we may view $\BRA
{\xi}$ as an operator $\E \to L^2 (G, B)$.  The adjoint of $\BRA {\xi}$ exists
and extends $\KET {\xi}$ to an operator $L^2 (G, B) \to \E$.  Let $\AdjG (L^2
(G, B), \E)$ be the space of equivariant, adjointable operators $L^2 (G, B)
\to \E$.  Then
$$
\E \supseteq \E_\si \cong \KET {\E_\si} \subseteq \AdjG (L^2 (G, B), \E).
$$
If~$G$ is compact, then $\E = \E_\si$.  If~$G$ is discrete, then $\KET
{\E_\si} = \AdjG (L^2 (G, B), \E)$.  We examine the operators $\KET {\xi}$ and
$\BRA {\xi}$ and the space $\E_\si$ in detail in Section~\ref{sec:si}.

We call~$\E$ \emph{square-integrable} iff $\E_\si$ is dense in~$\E$.  We
call~$B$ \emph{square-integrable} iff it is square-integrable as a Hilbert
module over itself.  Square-integrable \Cstar{}algebras are called ``proper''
in \cite{rieffel:si} and~\cite{meyer:KKG}.  The name ``square-integrable'' is
motivated by the relationship to square-integrable group representations
observed by Rieffel~\cite{rieffel:si}.  Square-integrable Hilbert modules are
characterized by the existence of many equivariant, adjointable operators to
$L^2 (G, B)$.  This gives rise to an equivariant version of Kasparov's
Stabilization Theorem~\cite{meyer:KKG}: A countably generated Hilbert module
is a direct summand of $L^2 (G, B)^\infty$ if and only if it is
square-integrable.

The fundamental example of a square-integrable Hilbert module is $L^2 (G, B)$.
All elements of $C_c (G, B)$ are square-integrable.  The closure of
$$
\KET {C_c (G, B)} \defeq \{ \KET {K} \mid K \in C_c (G, B) \}
$$
may be identified with the reduced crossed product $\Cred (G, B)$.  We always
think of $\Cred (G, B)$ as a subalgebra of $\AdjG \bigl (L^2 (G, B) \bigr)$ in
this way.

Our notation emphasizes that $\KET {\xi}$ is part of an inner product $\BRAKET
{\xi} {\eta} \defeq \BRA {\xi} \circ \KET {\eta}$ on $\E_\si$.  Since we want
Hilbert modules over $\Cred (G, B)$, we need subsets $\RC \subseteq \E_\si$
for which $\BRAKET {\RC} {\RC}$ is contained in $\Cred (G, B)$.  Following
Exel~\cite{exel:spectral}, we call such a subset \emph{relatively continuous}.
Let $\F = \F (\E, \RC) \subseteq \AdjG (L^2 (G, B), \E)$ be the closed linear
span of $\KET {\RC} \circ \Cred (G, B)$.  Then
\begin{equation}  \label{eq:concrete_i}
  \F \circ \Cred (G, B) \subseteq \F, \qquad
  \F^\ast \circ \F \subseteq \Cred (G, B),
\end{equation}
so that~$\F$ becomes a Hilbert module over $\Cred (G, B)$ with respect to the
inner product $\braket {\xi} {\eta} \defeq \xi^\ast \circ \eta$ and right
module structure $\xi \cdot x \defeq \xi \circ x$.  The closed linear span
$\Fix (\E, \RC)$ of $\F \circ \F^\ast \subseteq \AdjG (\E)$ is the
\emph{generalized fixed point algebra}.

There is a canonical isomorphism between $\Fix (\E, \RC)$ and the
\Cstar{}algebra $\Comp (\F)$ of compact operators on~$\F$.  Thus $\Fix (\E,
\RC)$ is Morita-Rieffel equivalent to an ideal in $\Cred (G, B)$.  For
compact~$G$, we get an ordinary fixed point algebra because $\Fix (\E,\E) =
\Comp (\E)^G$.

To exclude degenerate cases, we usually assume that $\RC$ is dense in~$\E$.

There are two obvious questions: Is square-integrability enough to guarantee
the existence of a dense, relatively continuous subspace $\RC \subseteq \E$?
Are $\F (\E, \RC)$ and $\Fix (\E, \RC)$ independent of~$\RC$?  Unfortunately,
the answer to both questions is negative.  Counterexamples come from
square-integrable representations of Abelian groups on Hilbert space.  This
situation can be analyzed completely (Section~\ref{sec:counter}).  The
problems are due to the subtle difference between continuous and measurable
fields of Hilbert spaces over the dual group.

Positive answers can be obtained if we require much more than
square\Mp{{}}integrability.  Following Kasparov, we call~$B$ \emph{proper} iff
there are a proper \Mpn{G}space~$X$ and an essential, equivariant
\Mpn{\ast}homomorphism from $C_0 (X)$ into the \emph{center} of the multiplier
algebra $\Mult (B)$ of~$B$.  If~$B$ is proper, then any Hilbert
\Mp{B,G}module~$\E$ is square-integrable, and $\F (\E, \RC)$ and $\Fix (\E,
\RC)$ do not depend on~$\RC$.  Actually, we can do with slightly less than
properness.  We only need that the induced action of~$G$ on the primitive
ideal space of~$B$ is proper.

Our main theoretical result is that the construction $(\E, \RC) \mapsto \F
(\E, \RC)$ can be inverted.  That is, all Hilbert modules over $\Cred (G, B)$
are of the form $\F (\E, \RC)$ for suitable $(\E, \RC)$, and $(\E, \RC)$ is
unique up to isomorphism if we impose further conditions on~$\RC$.  Let~$\F$
be a Hilbert module over $\Cred (G, B)$.  Define
\begin{equation}  \label{eq:EF}
  \E \defeq \F \otimes_{\Cred (G, B)} L^2 (G, B),
  \qquad
  \RC \defeq \F \otimes C_c (G, B) \subseteq \E.
\end{equation}
The \Mp{G}action on~$\E$ comes from the trivial action on~$\F$ and the usual
action on $L^2 (G, B)$.  Then~$\RC$ is dense in~$\E$ and relatively
continuous, and $\F (\E, \RC) \cong \F$.

The basis for our work is a detailed analysis of the construction of $\F (\E,
\RC)$.  It splits into two parts.  First, a relatively continuous subset $\RC
\subseteq \E$ yields a closed linear subspace $\F \subseteq \AdjG (L^2 (G, B),
\E)$ satisfying~\eqref{eq:concrete_i}.  The key idea here is the map $\xi
\mapsto \KET {\xi}$.  Secondly, $\F$ is turned into a Hilbert module over
$\Cred (G, B)$.  Only the first part uses special properties of groups.  The
second part should work equally well for coactions or actions of Hopf
algebras.  To simplify future extensions of this kind, we treat the second
part in greater generality.

Namely, we replace $L^2 (G, B)$ and $\Cred (G, B)$ by $\EL$ and~$A$,
where~$\EL$ is any Hilbert \Mp{B,G}module and $A$ is an essential
\Cstar{}subalgebra of $\AdjG (\EL)$.  A closed linear subspace $\F \subseteq
\AdjG (\EL, \E)$ is called a \emph{concrete Hilbert \Mpn{A}module} iff $\F
\circ A \subseteq \F$ and $\F^\ast \circ \F \subseteq A$.  We call~$\F$
\emph{essential} iff the linear span of $\F (\EL)$ is dense in~$\E$.  A
concrete Hilbert \Mpn{A}module carries a canonical Hilbert \Mpn{A}module
structure.  We view the embedding $\F \subseteq \AdjG (\EL, \E)$ as a
\emph{representation} of~$\F$.  The definition of a concrete Hilbert module is
relative to the representation $A \subseteq \AdjG (\EL)$.  This has the
consequence that all \emph{essential} representations of~$\F$ are isomorphic
to a canonical representation
$$
\F \cong \Comp ( A, \F) \subseteq
\AdjG (A \otimes_A \EL, \F \otimes_A \EL) \cong
\AdjG (\EL, \F \otimes_A \EL).
$$
In particular, $\E \cong \F \otimes_A \EL$ if $\F \subseteq \AdjG (\EL, \E)$
is an essential, concrete Hilbert module over~$A$.  The space $\F (\E, \RC)
\subseteq \AdjG (\EL, \E)$ is an essential, concrete Hilbert module over
$\Cred (G, B)$ if $\RC \subseteq \E$ is dense and relatively continuous.  This
explains the first half of~\eqref{eq:EF}.

Let~$\F$ be a concrete Hilbert module over $\Cred (G, B)$ and let $\RC_\F$ be
the set of all $\xi \in \E_\si$ with $\KET {\xi} \in \F$.  Then $\F (\E,
\RC_\F) = \F$.  A subset $\RC \subseteq \E$ is of the form $\RC_\F$ for some
concrete Hilbert module~$\F$ over $\Cred (G, B)$ if and only if it is
relatively continuous and \emph{complete} in an appropriate sense.  In
addition, $\F$ is essential iff $\RC_\F$ is dense.  A Hilbert module~$\E$
equipped with a dense, complete, relatively continuous subspace $\RC \subseteq
\E$ is called a \emph{continuously square-integrable Hilbert module}.  This
name is motivated by the case $B = \C$ and~$G$ Abelian, where~$\RC$ allows to
recover a continuous field of Hilbert spaces from a measurable field.  Our
analysis shows that $(\E, \RC) \mapsto \F (\E, \RC)$ yields a bijection
between isomorphism classes of continuously square-integrable Hilbert
\Mp{B,G}modules and isomorphism classes of Hilbert modules over $\Cred (G,
B)$.

For trivial coefficients and groupoids instead of groups, the correspondence
between Hilbert modules over $\Cred G$ and continuously square-integrable
representations of~$G$ on Hilbert space has been observed already by
Connes~\cite{connes:foliations}.  In order to do index theory on foliated
manifolds, he has to deal with the reduced \Cstar{}algebra of the holonomy
groupoid~$G$ of the foliation and the \Mp{KK}theoretic description of its
\Mpn{K}theory as $K_\ast (\Cred G) \cong KK_\ast (\C, \Cred G)$.  However, we
know very little about Hilbert modules over $\Cred G$.  Already the
determination of $K_\ast (\Cred G)$ is a major problem.  Therefore, Connes
replaces a Hilbert module over $\Cred G$ by a square-integrable representation
on a Hilbert space, equipped with a dense subspace with suitable properties.

Although most of our results can be extended to groupoids \textit{mutatis
mutandis}, we have decided not to cover groupoids.  Otherwise, we would have
to translate the basic theory of square-integrability, which so far has been
written down only for groups, and to use more complicated notation, since
groupoids do not act on \Cstar{}algebras but on bundles of \Cstar{}algebras.
These changes would make the article more difficult to read without changing
the content significantly.

A basic observation about square-integrable Hilbert modules is that~$\E$ is
square-integrable iff $\Comp (\E)$ is square-integrable.  This continues to
hold for continuously square-integrable Hilbert modules.  We can construct
relatively continuous subspaces of $\Comp (\E)$ from relatively continuous
subspaces of~$\E$ and vice versa.  These constructions are inverse to each
other if the group~$G$ is \emph{exact}, that is, the functor $\Cred (G,
\blank)$ preserves short exact sequences.  Otherwise, not all relatively
continuous subspaces of~$\E$ come from $\Comp (\E)$.  Since Abelian groups are
exact, the counterexamples in Section~\ref{sec:counter}---which involve group
actions on Hilbert spaces to begin with---also yield counterexamples in the
realm of group actions on \Cstar{}algebras.

\section{Some conventions}
\label{sec:conventions}

Throughout this article, $G$ is a locally compact group, $B$ is a
\Cstar{}algebra equipped with a strongly continuous action $\beta \colon G \to
\Aut (B)$ of~$G$ or, briefly, a \emph{\Mpn{G}\Cstar{}algebra}, and~$\E$ is a
\Mpn{G}equivariant Hilbert module over~$B$ with \Mpn{G}action $\gamma$ or,
briefly, a \emph{Hilbert \Mp{B,G}module}.  See, for instance, \cite{meyer:KKG}
for details.  We denote elements of~$G$ by $g$, $g'$, $g_2$, and fix a left
invariant \emph{Haar measure} $dg$ on~$G$.  Let $L^2G \defeq L^2 (G, dg)$ and
let~$G$ act on $L^2G$ via the left regular representation.  Let $\Delta \colon
G \to \R^\ast_+$ be the \emph{modular function} of~$G$ with the conventions
$d(g^{-1}) = \Delta (g^{-1}) \,dg$ and $d(gg_2) = \Delta (g_2) \,dg$.  We
write $\Adj (\E)$ and $\Comp (\E)$ for the \Cstar{}algebras of
\emph{adjointable} and \emph{compact} operators on~$\E$ and denote the
\Cstar{}algebra of \Mpn{G}equivariant, adjointable operators on~$\E$ by $\AdjG
(\E)$.

Since \emph{tensor products of Hilbert modules} are very important for us, we
recall the definition.  Let~$A$ be another \Mpn{G}\Cstar{}algebra, let~$\E_1$
be a Hilbert \Mpn{A,G}module and let~$\E$ be a Hilbert \Mpn{B,G}module.  Let
$\phi \colon A \to \Adj (\E)$ be an equivariant \Mpn{\ast}homomorphism.  Then
$\E_1 \otimes_A \E = \E_1 \otimes_\phi \E$ is a Hilbert \Mpn{B,G}module.  It
is the completion of the algebraic tensor product $\E_1 \otimes^\alg \E$ with
respect to the inner product
\begin{equation}  \label{eq:def_tensor}
  \braket {x_1 \otimes \xi_1} {x_2 \otimes \xi_2} \defeq
  \braket {\xi_1} {\phi(\braket {x_1} {x_2}_A) \xi_2}
  \qquad \forall\ x_1, x_2 \in \E_1,\ \xi_1, \xi_2 \in \E.
\end{equation}
The group~$G$ acts diagonally on $\E_1 \otimes_A \E$ by $\gamma_g (x \otimes
\xi) \defeq \gamma^{\E_1}_g(x) \otimes \gamma_g (\xi)$.  If $A = \C$, we drop
it from our notation.  For $\E_1 = L^2G$, we get the Hilbert \Mpn{B,G}module
$$
L^2 (G, B) \defeq L^2 G \otimes B.
$$

The \emph{bra-ket notation} is very useful in connection with Hilbert modules.
For $\xi \in \E$ we define the operators $\ket {\xi} \colon B \to \E$ and
$\bra {\xi} \colon \E \to B$ by $\ket {\xi} (b) \defeq \xi \cdot b$ and $\bra
{\xi} (\eta) \defeq \braket {\xi} {\eta}$, respectively.  These operators are
adjoints of one another.  The composition $\bra {\xi} \circ \ket {\eta}$ is
the operator of multiplication with the inner product $\braket {\xi} {\eta}$.
The composition $\ket {\xi} \circ \bra {\eta}$ is the ``rank-one operator''
$\ket {\xi} \bra {\eta} (\zeta) \defeq \xi \cdot \braket {\eta} {\zeta}$.

The map $\xi \mapsto \ket {\xi}$ is an isomorphism from~$\E$ onto $\Comp (B,
\E)$.  The map $\xi \mapsto \bra {\xi}$ is an isomorphism from the dual
$\E^\ast$ of~$\E$, which is a Hilbert module over $\Comp (\E)$ with $\Comp
(\E^\ast) \cong B$, onto $\Comp (\E, B)$.  The map
\begin{equation}  \label{eq:Comp_ketbra}
  \E \otimes_B \E^\ast \to \Comp (\E),
  \qquad \xi \otimes \eta \mapsto \ket {\xi} \bra {\eta},
\end{equation}
is an isomorphism of Hilbert modules over $\Comp (\E)$.

\section{The reduced crossed product}
\label{sec:crossed}

One of the basic observations of~\cite{rieffel:proper} is that $\Cred (G, B)$
arises as the generalized fixed point algebra of $\Comp (L^2 G) \otimes B$.
In order to make this isomorphism straightforward, we leave out the modular
function in the adjoint.  With this convention, modular functions do not show
up in any of our formulas.  Furthermore, it becomes easier to extend our
results to groupoids, where our convention is the standard one.  In order to
help readers who prefer the other convention, we explain how to modify
formulas if the adjoint is defined differently.

Let $C_c (G, B)$ be the space of continuous functions from~$G$ to~$B$ with
compact support.  We equip $C_c (G, B)$ with the following $\ast$-algebra
structure:
\begin{align}
  \label{eq:Cred_product}
  K \ast L (g) &\defeq
  \int_G K (g_2) \cdot \beta_{g_2} \bigl( L(g_2^{-1} g) \bigr) \,dg_2,
  \\
  \label{eq:Cred_adj}
  K^\ast (g) &\defeq
  \beta_g \bigl( K(g^{-1}) \bigr)^\ast.
  \\
\intertext{Usually, the adjoint is defined by}
  \label{eq:Cred_adj_mod}
  K^\times (g) &\defeq
  \beta_g \bigl( K(g^{-1}) \bigr)^\ast \cdot \Delta (g)^{-1}.
\end{align}
Equations \eqref{eq:Cred_adj} and~\eqref{eq:Cred_adj_mod} yield isomorphic
$\ast$-algebras.  The isomorphism is the map
$$
\mu \colon C_c (G, B) \to C_c (G, B),
\qquad (\mu K) (g) \defeq \Delta (g)^{1/2} K(g).
$$
Straightforward computations show
$$
\mu (K \ast L) = \mu (K) \ast \mu (L),
\qquad
(\mu K)^\ast = \mu (K^\times)
\qquad \forall\ K, L \in C_c (G, B).
$$

The advantage of~\eqref{eq:Cred_adj} is that the formula
\begin{equation}  \label{eq:Cred_act}
  (\rho_K f) (g) \defeq
  \int_G \beta_g \bigl( K(g^{-1} g_2) \bigr) \cdot f (g_2) \,dg_2
  \qquad \forall\ g \in G,\ K, f \in C_c (G, B)
\end{equation}
defines a \Mpn{\ast}homomorphism $\rho \colon C_c (G, B) \to \AdjG \bigl( L^2
(G, B) \bigr)$.  If we define the adjoint by~\eqref{eq:Cred_adj_mod}, we must
replace~$\rho$ by $\rho \circ \mu$.  That is, we have to insert a factor of
$\Delta (g^{-1} g_2)^{1/2}$ under the integral in~\eqref{eq:Cred_act}.

The \emph{reduced crossed product} $\Cred (G, B)$ is defined as the closure of
$\rho \bigl( C_c (G, B) \bigr)$ with respect to the operator norm on $\AdjG
\bigl( L^2 (G, B) \bigr)$.

If an adjointable operator on $L^2 (G, B)$ satisfies~\eqref{eq:Cred_act} for
some not necessarily compactly supported continuous function $K \colon G \to
B$, then we call it a \emph{Laurent operator} with \emph{symbol}~$K$
(following Exel's notation for Abelian groups~\cite{exel:spectral}).  If we
define the adjoint by~\eqref{eq:Cred_adj_mod}, then~$\rho$ is replaced by
$\rho \circ \mu$.  As a result, symbols are multiplied pointwise by $\Delta
(g)^{-1/2}$.

We can also define~$\rho_K$ if~$K$ is only a distribution on~$G$ taking values
in $\Mult (B)$.  In particular, we consider the distributions $b \cdot
\delta_1$ for $b \in B$ and $\delta_g$ for $g \in G$ that are defined by
$\int_G b \cdot \delta_1 (g) \cdot f(g) \,dg \defeq b f(1)$ and $\int_G
\delta_g (g_2) \cdot f(g_2) \,dg_2 = f(g)$.  If we plug them
into~\eqref{eq:Cred_act}, we get the operators $\rho_b, \rho_g \in \AdjG
\bigl( L^2 (G, B) \bigr)$,
\begin{align}
  \label{eq:rho_B}
  \rho_b (f) (g_2) &\defeq \beta_{g_2} (b) \cdot f(g_2),
  \\
  \label{eq:rho_G}
  \rho_g (f) (g_2) &\defeq f(g_2 g).
\end{align}
We have $\rho_g^\ast = \rho_{g^{-1}} \cdot \Delta (g)^{-1}$, so that $\rho_g
\cdot \Delta (g)^{1/2}$ is unitary.  It is elementary to verify that $\rho_b$
and~$\rho_g$ multiply $\rho \bigl( C_c (G, B) \bigr)$ and hence are contained
in $\Mult \bigl( \Cred (G, B) \bigr)$.  If we define the adjoint
by~\eqref{eq:Cred_adj_mod} and replace~$\rho$ by $\rho \circ \mu$,
then~$\rho_b$ remains unchanged and~$\rho_g$ is replaced by $\rho_g \cdot
\Delta (g)^{1/2}$.

We can view $C_c (G, B)$ as the inductive limit of the Banach spaces of
continuous functions $G \to B$ with support contained in $K \subseteq G$,
where~$K$ runs through the compact subsets of~$G$.  Hence $C_c (G, B)$ is a
complete bornological vector space and thus a complete topological vector
space in a canonical way.

The usual \Mp{L^1}norm on $C_c (G, B)$ is \emph{not} a \Mpn{\ast}algebra norm
and does not control the operator norm on $\Cred (G, B)$ because we left out
the modular function in the adjoint.  As a substitute, we define
$$
\norm {K}_I \defeq
\max \int_G \norm {K (g)} \,dg,\ \int_G \norm {K^\ast (g)} \,dg.
$$
This norm is submultiplicative and satisfies $\norm {K^\ast}_I = \norm {K}_I$
by definition.  An application of the Cauchy-Schwarz inequality shows that
\begin{multline*}
  \norm {\mu^{-1}(K)}_1
  \defeq
  \int_G \norm {\mu^{-1} (K) (g)} \,dg
  =
  \int_G \norm {K (g)} \cdot \Delta (g)^{-1/2} \,dg
  \\ \le
  \left( \int_G \norm {K (g)} \,dg \right)^{1/2} \cdot
  \left( \int_G \norm {K (g)} \,d (g^{-1}) \right)^{1/2}
  \le
  \norm {K}_I.
\end{multline*}

Hence $\norm {\rho_K} \le \norm {K}_I$.  The latter estimate continues to
hold for groupoids~\cite{renault:representations}.  Finally, we mention that
$C_c (G, B)$ has approximate identities:

\begin{lemma}  \label{lem:apprid}
  There is a net $(u_j)_{j \in J}$ of elements of $C_c (G, B)$ such that:
  \begin{itemize}%

  \item $u_j = u_j^\ast$ for all $j \in J$;

  \item $(u_j)$ is bounded with respect to the norm $\norm {\blank}_I$;

  \item $(u_j)$ is an approximate identity of $C_c (G, B)$ and of $\Cred (G,
    B)$ with respect to the inductive limit bornology and the operator norm,
    respectively.

  \end{itemize}
\end{lemma}

\begin{proof}
  This follows from Lemma~3.2 of~\cite{renault:representations}.
\end{proof}

\section{Square-integrable Hilbert modules}
\label{sec:si}

In the introduction, we called $\xi \in \E$ square-integrable iff $\BRA {\xi}
\eta \in L^2 (G, B)$ for all $\eta \in \E$.  We have to explain what $\BRA
{\xi} \eta \in L^2 (G, B)$ means.  Let $(\chi_i)_{i \in I}$ be a net of
continuous, compactly supported functions $G \to [0,1]$ with $\chi_i (g) \to
1$ uniformly on compact subsets of~$G$.  Let $f \in C_b (G, B)$.  We say
that~$f$ is \emph{square-integrable} and write $f \in L^2 (G, B)$ iff the net
$(\chi_i \cdot f)_{i \in I}$ converges in $L^2 (G, B)$.  We identify~$f$ with
the limit of this net, so that~$f$ becomes an element of $L^2 (G, B)$.

As a result, we may view $\BRA {\xi}$ as an operator $\E \to L^2 (G, B)$ if
$\xi \in \E_\si$.  The closed graph theorem implies that $\BRA {\xi}$ is
bounded as a map to $L^2 (G, B)$.  Since $\BRA {\xi}$ is bounded as an
operator to $C_b (G, B)$, its graph in $\E \times L^2 (G, B)$ is closed.
(Since we do not assume~$G$ to be \Mpn{\sigma}compact as in~\cite{meyer:KKG},
we cannot employ the Banach-Steinhaus theorem as in the proof of \cite
[Lemma~8.1]{meyer:KKG}.)  The same argument as in the proof of \cite[Lemma
8.1]{meyer:KKG} shows that $\BRA {\xi}$ is adjointable and that its adjoint
extends $\KET {\xi}$ to an operator $L^2 (G, B) \to \E$.  Conversely, suppose
that the operator $\KET {\xi}$ defined in~\eqref{eq:def_KET} extends to an
\emph{adjointable} operator $L^2 (G, B) \to \E$.  The computation that yields
$\KET {\xi} = \BRA {\xi}^\ast$ shows that $\KET {\xi}^\ast (\eta) = \BRA {\xi}
\eta \in L^2 (G, B)$ for all $\eta \in \E$.  Hence $\xi \in \E_\si$ iff $\KET
{\xi}$ extends to an adjointable operator $L^2 (G, B) \to \E$.

It is clear that $\E_\si$ is a vector space.  It becomes a Banach space when
equipped with the norm
$$
\norm{\xi}_\si \defeq
\norm {\braket {\xi} {\xi}}^{1/2} + \norm {\KET {\xi}} =
\norm {\braket {\xi} {\xi}}^{1/2} + \norm {\BRAKET {\xi} {\xi}}^{1/2}.
$$

The remainder of this section contains elementary computations with the
operators $\KET {\xi}$ and $\BRA {\xi}$ that are needed later.  It is
convenient for reference purposes to collect these computation in a single
section.  We have
\begin{xalignat}{2}
  \label{eq:KET_compose}
  \KET {T (\xi)} &= T \circ \KET {\xi}
  \qquad &\forall\ & T \in \AdjG (\E, \E'),\ \xi \in \E,
  \\
  \label{eq:KET_B}
  \KET {\xi \cdot b} &= \KET {\xi} \circ \rho_b
  \qquad &\forall\ & b \in B,\ \xi \in \E,
  \\
  \label{eq:KET_G}
  \KET {\gamma_g(\xi)} &= \KET {\xi} \circ \rho_g^\ast
  =
  \KET {\xi} \circ \rho_{g^{-1}} \cdot \Delta (g)^{-1}
  \qquad &\forall\ & g \in G,\ \xi \in \E.
\end{xalignat}
Equations \eqref{eq:KET_compose}--\eqref{eq:KET_G} follow at once from the
definitions \eqref{eq:def_KET}, \eqref{eq:rho_B}, and~\eqref{eq:rho_G}.
If~$\rho_g$ is replaced by $\rho_g \cdot \Delta (g)^{-1/2}$,
then~\eqref{eq:KET_G} has to be modified accordingly.  Therefore,
\begin{xalignat}{2}
  \label{eq:KET_compose_est}
  \norm {T (\xi)}_\si &\le \norm {\xi}_\si \cdot \norm {T}
  \qquad &\forall\ & T \in \AdjG (\E, \E'),\ \xi \in \E_\si,
  \\
  \label{eq:KET_B_est}
  \norm {\xi \cdot b}_\si &\le \norm {\xi}_\si \cdot \norm {b}
  \qquad &\forall\ & b \in B,\ \xi \in \E_\si,
  \\
  \label{eq:KET_G_est}
  \norm {\gamma_g (\xi)}_\si &\le
  \norm {\xi}_\si \cdot \max \{1, \Delta (g)^{-1/2} \}
  \qquad &\forall\ & g \in G,\ \xi \in \E_\si.
\end{xalignat}
Thus $\E_\si$ is \Mpn{G}invariant and a Banach bimodule over $\AdjG (\E)
\times B$.  However, the action of~$G$ on~$\E_\si$ need not be continuous, and
$\E_\si \cdot B$ need not be dense in $\E_\si$.

Let $\xi, \eta \in \E_\si$.  We compute the compositions $\KET {\xi} \BRA
{\eta}$ and $\BRAKET {\xi} {\eta}$.  Formally, we have $\KET {\xi} \BRA {\eta}
\zeta = \int_G \gamma_g (\xi) \braket {\gamma_g(\eta)} {\zeta} \,dg$.  To
interpret this integral, recall that $\braket {\gamma_g (\eta)} {\zeta}$ is
the limit of the net $(\chi_i (g) \cdot \braket {\gamma_g (\eta)} {\zeta})_{i
\in I}$ in $L^2 (G, B)$.  Hence
\begin{equation}  \label{eq:KETBRA}
  \KET {\xi} \BRA {\eta} =
  \int_G \gamma_g (\ket {\xi} \bra {\eta}) \,dg \defeq
  \lim_{i \in I} \int_G \chi_i (g) \cdot
  \gamma_g (\ket {\xi} \bra {\eta}) \,dg
\end{equation}
for all $\xi, \eta \in \E_\si$.  The limit exists in the strict
topology~\cite{meyer:KKG}.  Moreover,
$$
(\BRAKET {\xi} {\eta} f) (g) =
\braket {\gamma_g (\xi)} {\KET {\eta} f} =
\int_G \braket {\gamma_g (\xi)} {\gamma_{g_2} (\eta)} \cdot f(g_2)
\,dg_2
$$
for all $f \in C_c (G, B)$, $g \in G$.  Comparing with~\eqref{eq:Cred_act}, we
see that $\BRAKET {\xi} {\eta}$ is a Laurent operator, whose symbol $\BRAKET
{\xi} {\eta} \in C_b (G, B)$ is
\begin{equation}  \label{eq:BRAKET}
  \BRAKET {\xi} {\eta} (g) = \braket {\xi} {\gamma_g (\eta)}
  \qquad \forall\ g \in G,\ \xi,\eta \in \E_\si.
\end{equation}
It may happen that $\BRAKET {\xi} {\eta} \notin \Cred (G, B)$.  If we define
the adjoint by~\eqref{eq:Cred_adj_mod}, then~\eqref{eq:BRAKET} has to be
replaced by $\BRAKET {\xi} {\eta} (g) = \braket {\xi} {\gamma_g (\eta)} \cdot
\Delta (g)^{-1/2}$.

The basic example of a square-integrable Hilbert module is $L^2 (G, B)$.  We
claim that $C_c (G, B) \subseteq L^2 (G, B)_\si$.  Let $K \in C_c (G, B)$.
Then
$$
(\KET {K} f) (g) =
\int_G \bigl( \beta_{g_2} (K) \bigr) (g) \cdot f(g_2) \,dg_2 =
\int_G \beta_{g_2} \bigl( K(g_2^{-1} g) \bigr) \cdot f(g_2) \,dg_2.
$$
Comparing with \eqref{eq:Cred_act}, we see that
\begin{equation}  \label{eq:KET_check}
  \KET {K} = \rho_{\check{K}}, \qquad \KET {\check{K}} = \rho_K,
\end{equation}
if
\begin{equation}  \label{eq:def_check}
  \check {K} (g) \defeq \beta_g \bigl( K(g^{-1}) \bigr).
\end{equation}
The map $K \mapsto \check {K}$ is a bijection from $C_c (G, B)$ onto $C_c (G,
B)$.  As a result, $\KET {K}$ extends to an adjointable operator, so that $C_c
(G, B) \subseteq L^2 (G, B)_\si$ as asserted.

If we define the adjoint by~\eqref{eq:Cred_adj_mod}, then~$\rho$ has to be
replaced by $\rho \circ \mu$.  Hence we desire an equation $\KET {K} = \rho
\circ \mu (\tilde{K})$ instead of~\eqref{eq:KET_check} and put
\begin{equation}  \label{eq:def_tilde}
  \tilde{K} (g) \defeq
  \check{K} (g) \cdot \Delta (g)^{-1/2}
  = \beta_g \bigl( K(g^{-1}) \bigr) \cdot \Delta (g)^{-1/2}.
\end{equation}

We turn~$\E$ into a right module over the convolution algebra $C_c (G, B)$ by
\begin{equation}  \label{eq:right_module}
  \xi \ast K \defeq
  \KET {\xi} (\check{K}) =
  \int_G \gamma_g (\xi) \cdot \check{K} (g) \,dg =
  \int_G \gamma_g \bigl( \xi \cdot K(g^{-1}) \bigr) \,dg
\end{equation}
for all $\xi \in \E$, $K \in C_c (G, B)$.  Since $\KET {\xi}$ is equivariant,
\eqref{eq:KET_compose} and~\eqref{eq:KET_check} yield
\begin{equation}  \label{eq:KET_right_module}
  \KET {\xi \ast K} =
  \KET {\xi} \circ \KET {\check{K}} =
  \KET {\xi} \circ \rho_K
  \qquad \forall\ \xi \in \E,\ K \in C_c (G, B).
\end{equation}
Hence
$$
\KET {(\xi \ast K) \ast L} =
\KET {\xi} \circ \rho_K \circ \rho_L =
\KET {\xi} \circ \rho_{K \ast L} =
\KET {\xi \ast (K \ast L)}.
$$
Since the map $\xi \mapsto \KET {\xi}$ is injective, $\E$ is a right module
over $C_c (G, B)$.

If we define the adjoint by~\eqref{eq:Cred_adj_mod}, then we replace
$\check{K}$ by~$\tilde{K}$ in~\eqref{eq:right_module}.  The same computation
as above yields $\KET {\xi \ast K} = \KET {\xi} \circ \rho \mu(K)$ instead
of~\eqref{eq:KET_right_module}.

Using $\norm {\rho_K} \le \norm{K}_I$ and~\eqref{eq:KET_right_module}, we
obtain the following norm estimates:
\begin{align}
  \label{eq:right_estimate}
  \norm {\xi \ast K} &\le \norm {\xi} \cdot \norm {K}_I,
  \\
  \label{eq:right_estimate_si}
  \norm {\xi \ast K}_\si &\le \norm {\xi}_\si \cdot \norm {K}_I,
  \\
  \label{eq:right_estimate_si_KET}
  \norm {\xi \ast K}_\si &\le
  \norm {\KET {\xi}} \cdot
  \max \{ \norm {\rho_K}, \norm {\check{K}}_{L^2 (G, B)}  \}.
\end{align}

Lemma~8.1 of~\cite{meyer:KKG} asserts that~$\xi$ is contained in the closure
of $\KET {\xi} \bigl( C_c (G, B) \bigr)$.  Hence $\E \ast C_c (G, B)$ is dense
in~$\E$.  Using~\eqref{eq:right_estimate}, we conclude that the approximate
identities $(u_j)_{j \in J}$ of Lemma~\ref{lem:apprid} satisfy
\begin{equation}  \label{eq:apprid_right_module}
  \lim {} \norm {\xi \ast u_j - \xi} = 0
  \qquad
  \forall\ \xi \in \E.
\end{equation}
However, $\E_\si \ast C_c (G, B)$ need not be dense in $\E_\si$ with respect
to $\norm {\blank}_\si$.

\section{Representations of Hilbert modules}
\label{sec:concrete_Hilbert}

Throughout this section, we let~$\EL$ be a \Mpn{G}equivariant Hilbert module
over a \Mpn{G}\Cstar{}algebra~$B$, and we let $A \subseteq \AdjG (\EL)$ be an
\emph{essential} \Cstar{}subalgebra.  That is, the closed linear span of $A
\cdot \EL$ is dense in~$\EL$.  By Cohen's Factorization Theorem, this
implies $A \cdot \EL = \EL$.  We are particularly interested in the case
$A = \Cred (G, B)$, $\EL = L^2 (G, B)$.  The group~$G$ is only there because
we want to assert that our constructions are invariant with respect to a group
action.

\begin{definition}  \label{def:concrete_Hilbert}
  Let~$\E$ be a Hilbert \Mp{B,G}module.  A \emph{concrete Hilbert
  \Mpn{A}module} is a closed linear subspace $\F \subseteq \AdjG (\EL, \E)$
  that satisfies $\F \circ A \subseteq \F$ and $\F^\ast \circ \F \subseteq A$.
  We call~$\F$ \emph{essential} iff the linear span of $\F (\EL)$ is dense
  in~$\E$.
\end{definition}

A concrete Hilbert \Mpn{A}module $\F \subseteq \AdjG (\EL, \E)$ can be made
essential by making~$\E$ smaller.  Let $\E' \subseteq \E$ be the closed linear
span of $\F (\EL)$.  Then $\E'$ is an invariant Hilbert submodule and $\F
\subseteq \AdjG (\EL, \E')$ is an essential, concrete Hilbert \Mpn{A}module.

\begin{lemma}  \label{lem:concrete_Hilbert}
  Let $\F \subseteq \AdjG (\EL, \E)$ be a concrete Hilbert \Mpn{A}module.
  Then~$\F$ becomes a Hilbert \Mp{A}module when equipped with the right
  \Mpn{A}module structure
  $$
  \xi \cdot a \defeq \xi \circ a
  \qquad
  \forall\ \xi \in \F,\ a \in A
  $$
  and the \Mpn{A}valued inner product
  $$
  \braket {\xi} {\eta} \defeq \xi^\ast \circ \eta
  \qquad
  \forall\ \xi, \eta \in \F.
  $$
  The Hilbert module norm and the operator norm on~$\F$ coincide.  We have
  \begin{equation}  \label{eq:Fthree}
    \F = \F \circ A = \F \circ \F^\ast \circ \F
  \end{equation}
  and\/ $\F (\EL) = \F \circ \F^\ast (\E) = \F \circ \F^\ast \circ \F (\EL)$.
  Hence~$\F$ is essential iff the linear span of\/ $\F \circ \F^\ast (\E)$ is
  dense in~$\E$.
\end{lemma}

We always furnish a concrete Hilbert \Mpn{A}module with the Hilbert
\Mpn{A}module structure defined above.

\begin{proof}
  By assumption, $\xi \cdot a \in \F$ for all $\xi \in \F$, $a \in A$ and
  $\braket {\xi} {\eta} \in A$ for all $\xi, \eta \in \F$.  The conditions
  $$
  \braket {\xi} {\eta \cdot a} = \braket {\xi} {\eta} \cdot a,
  \quad
  \braket {\xi} {\eta} = (\braket {\eta} {\xi})^\ast,
  \quad
  \braket {\xi} {\xi} \ge 0
  $$
  for a pre-Hilbert module are obviously satisfied.  Since
  \begin{equation}  \label{eq:T_isometric}
    \norm {\xi} =
    \norm {\xi^\ast \xi}^{1/2} =
    \norm {\braket {\xi} {\xi}}^{1/2},
  \end{equation}
  the norm that comes from the \Mpn{A}valued inner product equals the operator
  norm.  Hence~$\F$ is a Hilbert module.  We have $\F \circ \F^\ast \circ \F
  \subseteq \F \circ A \subseteq \F$.  It is a general feature of Hilbert
  modules that any $\xi \in \F$ may be written as $\eta \braket {\eta} {\eta}$
  for some $\eta \in \F$ \cite[Lemme~1.3]{blanchard}.  Hence $\F \subseteq \F
  \circ \F^\ast \circ \F$.  Equation~\eqref{eq:Fthree} follows.  Since $\F
  (\EL) = \F \circ \F^\ast \circ \F (\EL) \subseteq \F \circ \F^\ast (\E)
  \subseteq \F (\EL)$, these three sets are equal.
\end{proof}

Let~$\F$ be a Hilbert \Mpn{A}module.  We construct a canonical representation
of~$\F$ as a concrete Hilbert module.  We equip~$\F$ with the trivial action
of~$G$, so that $\F \otimes_A \EL$ is a Hilbert \Mp{B,G}module.  The
construction $\F \mapsto \F \otimes_A \EL$ is functorial, that is, an
adjointable operator $x \colon \F_1 \to \F_2$ between Hilbert \Mpn{A}modules
induces an equivariant, adjointable operator $x \otimes \ID_\EL \colon \F_1
\otimes_A \EL \to \F_2 \otimes_A \EL$.  Using the isomorphism $A \otimes_A \EL
\cong A \cdot \EL = \EL$, we obtain a map
\begin{equation}  \label{eq:def_T}
  T \colon \F \cong \Comp (A, \F) \longrightarrow
  \AdjG \bigl( A \otimes_A \EL, \F \otimes_A \EL \bigr) \cong
  \AdjG \bigl( \EL, \F \otimes_A \EL \bigr).
\end{equation}
More explicitly, we have $T (\xi) (f) \defeq \xi \otimes f$ and $T (\xi)^\ast
(\eta \otimes f) \defeq \braket {\xi} {\eta} (f)$ for all $\xi, \eta \in \F$
and $f \in \EL$, where we view $\braket {\xi} {\eta} \in A \subseteq \AdjG
(\EL)$.

\begin{theorem}  \label{the:concrete_Hilbert}
  Let~$\F$ be a Hilbert \Mpn{A}module and define~$T$ as in~\eqref{eq:def_T}.
  Then $T (\F)$ is an essential, concrete Hilbert \Mpn{A}module and $T \colon
  \F \to T (\F)$ is an isomorphism of Hilbert \Mpn{A}modules.  If\/ $\F
  \subseteq \AdjG (\EL, \E)$ already is an essential, concrete Hilbert
  \Mpn{A}module, then
  $$
  U \colon \F \otimes_A \EL \to \E,
  \qquad \xi \otimes f \mapsto \xi (f),
  $$
  is an equivariant unitary that satisfies $U \circ \bigl( T (\xi) \bigr) =
  \xi$ for all $\xi \in \F$.  That is, $\F$ and $T (\F)$ are isomorphic as
  concrete Hilbert \Mpn{A}modules via~$U$.
\end{theorem}

\begin{proof}
  We have $T (\xi \cdot a) = T (\xi) \circ a$ and $T (\xi)^\ast T (\eta) =
  \braket {\xi} {\eta}$ for all $\xi, \eta \in \F$, $a \in A$.
  Equation~\eqref{eq:T_isometric} shows that~$T$ is isometric, so that $T
  (\F)$ is closed.  Thus $T(\F)$ is a concrete Hilbert \Mpn{A}module and $T
  \colon \F \to T(\F)$ is an isomorphism with respect to the Hilbert
  \Mpn{A}module structure of Lemma~\ref{lem:concrete_Hilbert}.  $T(\F)$ is
  essential because $\F \otimes_A \EL$ is generated by elementary tensors $\xi
  \otimes f = T(\xi) (f)$ with $\xi \in \F$ and $f \in \EL$.

  Suppose that $\F \subseteq \AdjG (\EL, \E)$ is a concrete Hilbert
  \Mpn{A}module.  The map~$U$ is isometric (hence well-defined)
  by~\eqref{eq:def_tensor} and equivariant.  If~$\F$ is essential, then the
  range of~$U$ is dense, so that~$U$ is unitary.  We compute $U \bigl( T (\xi)
  (f) \bigr) = U (\xi \otimes f) = \xi (f)$ for all $\xi \in \F$, $f \in \EL$.
  That is, $U \circ \bigl( T (\xi) \bigr) = \xi$.
\end{proof}

Put in a nutshell, any Hilbert \Mpn{A}module~$\F$ can be represented as an
esssential, concrete Hilbert \Mpn{A}module, and this representation is unique
up to isomorphism.  The underlying Hilbert \Mp{B,G}module~$\E$ is canonically
isomorphic to $\F \otimes_A \EL$.

\begin{theorem}  \label{the:concrete_morphism}
  Let $\F \subseteq \AdjG (\EL, \E)$ be a concrete Hilbert \Mpn{A}module.  The
  map
  $$
  \ket {\xi} \bra {\eta} \mapsto \xi \circ \eta^\ast
  \in \F \circ \F^\ast \subseteq \AdjG (\E)
  $$
  extends to a \Mpn{\ast}isomorphism from $\Comp (\F)$ onto the norm closure
  of\/ $\F \circ \F^\ast$ in $\AdjG (\E)$.  This representation of\/ $\Comp
  (\F)$ is essential iff\/~$\F$ is essential.

  If\/~$\F$ is essential, we may extend this representation of\/ $\Comp (\F)$
  to a strictly continuous, injective, unital \Mpn{\ast}homomorphism $\phi
  \colon \Adj (\F) \to \AdjG (\E)$, whose range is
  $$
  M \defeq \{ x \in \Adj (\E) \mid
  x \circ \F \subseteq \F,\ x^\ast \circ \F \subseteq \F \}.
  $$
\end{theorem}

\begin{proof}
  It is clear that~$M$ is a \Cstar{}subalgebra of $\Adj (\E)$.  Let $D
  \subseteq \Adj (\E)$ be the closed linear span of $\F \circ \F^\ast$.  By
  construction, $D$ is closed and $D^\ast = D$.  Equation~\eqref{eq:Fthree}
  implies $D \circ \F \subseteq \F$ and hence $D \subseteq M$.  If $x \circ \F
  \subseteq \F$, then $x \circ D \subseteq D$.  Hence~$D$ is a closed ideal
  in~$M$.  Conversely, if $x \circ D \subseteq D$,
  then
  $$
  x \circ \F = x \circ \F \circ \F^\ast \circ \F \subseteq
  D \circ \F \subseteq \F
  $$
  by~\eqref{eq:Fthree}.  Consequently, $x \in M$ iff $xD \subseteq D$ and $Dx
  \subseteq D$.  We define a \Mpn{\ast}homomorphism $\psi \colon M \to \Adj
  (\F)$ by $\psi (x) (\xi) \defeq x \circ \xi$ for $x \in M$, $\xi \in \F$.
  If~$\F$ is essential, then~$\psi$ is injective because if $\psi (x) = 0$,
  then~$x$ vanishes on the dense subspace $\F (\EL) \subseteq \E$, so that $x
  = 0$.  In general, at least the restriction of~$\psi$ to~$D$ is injective
  because $x \circ \F = 0$ implies $x \circ D = 0$ and hence $xx^\ast = 0$.
  We have $\psi (\xi \circ \eta^\ast) = \ket {\xi} \bra {\eta}$ for all $\xi,
  \eta \in \F$.  Hence $\psi (D) = \Comp (\F)$ and~$\psi|_D^{-1}$ equals the
  map $\ket {\xi} \bra {\eta} \mapsto \xi \circ \eta^\ast$.
  Lemma~\ref{lem:concrete_Hilbert} implies that $\psi|_D^{-1} \colon
  \Comp (\F) \to \Adj (\E)$ is essential iff~$\F$ is essential.

  Therefore, if~$\F$ is essential, then $\psi|_D^{-1}$ extends uniquely to a
  strictly continuous, unital, injective \Mpn{\ast}homomorphism~$\phi$ from
  $\Adj (\F) \cong \Mult \bigl( \Comp (\F) \bigr)$ to $\Adj (\E)$.  The range
  of~$\phi$ is contained in~$M$ because $\Comp (\F)$ is an ideal in $\Adj
  (\F)$.  Since $\phi \circ \psi = \ID_M$, the map~$\phi$ is an isomorphism
  onto~$M$.
\end{proof}

If $\E = \F \otimes_A \EL$ and $\F \subseteq \AdjG (\EL, \E)$ is the standard
representation~\eqref{eq:def_T}, then the homomorphism $\phi \colon \Adj (\F)
\to \AdjG (\E)$ equals the canonical map $x \mapsto x \otimes_A \ID_\EL$.

\begin{corollary}  \label{cor:ideal}
  Let $\F \subseteq \AdjG (\EL, \E)$ be an essential, concrete Hilbert
  \Mpn{A}module.  The following statements are equivalent:
  \begin{enumerate}[(i)]%

  \item $\phi \colon \Adj (\F) \to \AdjG (\E)$ is an isomorphism onto $\AdjG
    (\E)$;

  \item we have $u \circ \F = \F$ for all $u \in \AdjG (\E)$;

  \item the closed linear span of\/ $\F \circ \F^\ast$ is an ideal in $\AdjG
    (\E)$.

  \end{enumerate}
\end{corollary}

We call~$\F$ \emph{ideal} iff one of these assertions holds.

\begin{proof}
  By Theorem~\ref{the:concrete_morphism}, condition~(i) is equivalent to $M =
  \AdjG (\E)$.  Condition~(ii) asserts that all unitaries $u \in \AdjG (\E)$
  are contained in~$M$.  Since any element of $\AdjG (\E)$ may be written as a
  sum of four unitaries, the first two conditions are equivalent.  In the
  proof of Theorem~\ref{the:concrete_morphism}, we observed that $x \in \Adj
  (\E)$ satisfies $xD \subseteq D$ and $Dx \subseteq D$ iff $x \in M$.  Hence
  condition~(iii) is equivalent to $M = \AdjG (\E)$ as well.
\end{proof}

\section{Continuously square-integrable Hilbert modules}
\label{sec:relcont}

It is convenient to keep the abbreviations
$$
A \defeq \Cred (G, B),
\qquad
\EL \defeq L^2 (G, B).
$$

\begin{definition}  \label{def:rel_cont}
  A subset $\RC \subseteq \E$ is called \emph{relatively continuous} iff $\RC
  \subseteq \E_\si$ and
  $$
  \BRAKET {\RC} {\RC} \defeq \{ \BRAKET {\xi} {\eta} \mid \xi, \eta \in \RC \}
  $$
  is contained in $\Cred (G, B) \subseteq \AdjG (\EL)$.  If $\RC \subseteq \E$
  is a relatively continuous subset, let $\F (\E, \RC) \subseteq \AdjG (\EL,
  \E)$ be the closed linear span of $\KET {\RC} \circ \Cred (G, B) \cup \KET
  {\RC}$.
\end{definition}

Recall that $\BRAKET {\xi} {\eta}$ is a Laurent operator, whose symbol is
given by~\eqref{eq:BRAKET}.  This often allows to verify relative continuity.
In many interesting examples we have $\norm { \BRAKET {\xi} {\eta} }_I <
\infty$ or even $\BRAKET {\xi} {\eta} \in C_c (G, B)$ for all $\xi, \eta \in
\RC$.

\begin{proposition}  \label{pro:rel_cont}
  Let $\RC \subseteq \E$ be relatively continuous.  Then $\F (\E, \RC)$ is a
  concrete Hilbert module over $\Cred (G, B)$.  If\/~$\RC$ is dense in~$\E$,
  then $\F (\E, \RC)$ is essential.
\end{proposition}

\begin{proof}
  By construction, $\F \defeq \F (\E, \RC)$ is a closed linear subspace and
  $\F \circ A \subseteq \F$.  The assumption $\BRAKET {\RC} {\RC} \subseteq A$
  implies $\F^\ast \circ \F \subseteq A$.  Suppose that~$\RC$ is dense
  in~$\E$.  Since $\E \ast C_c (G, B)$ is dense in~$\E$, the subset $\RC \ast
  C_c (G, B) = \KET {\RC} \bigl( C_c (G, B) \bigr)$ is dense in~$\E$.
  Therefore, $\F (\EL)$ is dense in~$\E$.
\end{proof}

\begin{proposition}  \label{pro:RC_dense}
  Let $\F \subseteq \AdjG (\EL, \E)$ be a concrete Hilbert module over $\Cred
  (G, B)$.  Define
  \begin{align*}
    \RC_\F &\defeq \{ x \in \E_\si \mid \KET {x} \in \F \},
    \\
    \RC_\F^0 &\defeq \{ \xi (K) \mid \xi \in \F,\ K \in C_c (G, B) \}.
  \end{align*}
  Then $\RC_F^0 \subseteq \RC_\F$.  Both $\RC_\F^0$ and $\RC_\F$ are
  relatively continuous, and\/ $\KET {\RC_\F^0}$ and\/ $\KET {\RC_\F}$ are
  dense in~$\F$.  Thus
  $$
  \F (\E, \RC_\F^0) = \F (\E, \RC_\F) = \F.
  $$
\end{proposition}

\begin{proof}
  It is evident that~$\RC_\F$ is relatively continuous.  Let $\xi \in \F$ and
  $K \in C_c (G, B)$.  Since~$\xi$ is equivariant, \eqref{eq:KET_compose}
  and~\eqref{eq:KET_check} yield
  $$
  \KET {\xi (K)} = \xi \circ \KET {K} = \xi \circ \rho_{\check{K}}
  \in \F \circ \Cred (G, B) \subseteq \F.
  $$
  This implies $\RC_\F^0 \subseteq \RC_\F$.  Thus $\RC_\F^0$ is relatively
  continuous.  The above computation shows $\KET {\RC_\F^0} = \F \cdot C_c (G,
  B)$.  Since $C_c (G, B)$ is dense in~$A$ and~$\F$ is a Hilbert
  \Mpn{A}module, $\F \cdot C_c (G, B)$ is dense in~$\F$.  It follows that
  $\KET {\RC_\F^0}$ and $\KET {\RC_\F}$ are dense subsets of~$\F$.  Therefore,
  $\F (\E, \RC_\F^0) = \F (\E, \RC_\F) = \F$.
\end{proof}

The subspace $\RC \subseteq \F \otimes_A \EL$ that is defined in~\eqref{eq:EF}
equals $\RC_\F^0$.

\begin{definition}  \label{def:complete}
  We call $\RC \subseteq \E$ \emph{complete} iff~$\RC$ is a linear subspace of
  $\E_\si$ that is closed with respect to the norm $\norm {\blank}_\si$ and
  satisfies $\RC \ast C_c (G, B) \subseteq \RC$.  The \emph{completion} of a
  subset $\RC \subseteq \E_\si$ is the smallest complete subset that
  contains~$\RC$.  That is, the completion of~$\RC$ is the \Mp{\norm
  {\blank}_\si}closed linear span of $\RC \cup \RC \ast C_c (G, B)$.

  A \emph{continuously square-integrable Hilbert \Mp{B,G}module} is a Hilbert
  \Mp{B,G}module together with a dense, complete, relatively continuous
  subspace.
\end{definition}

If $\RC \subseteq \E$ is a complete, relatively continuous subset, then the
closure of $\KET {\RC}$ is already a right \Mp{A}module
by~\eqref{eq:KET_right_module}.  Hence $\F (\E, \RC)$ is simply the closure of
$\KET {\RC}$.

\begin{theorem}  \label{the:RC_characterize}
  The map $\F \mapsto \RC_\F$ is a bijection from the set of concrete Hilbert
  \Mp{\Cred (G, B)}modules $\F \subseteq \AdjG (L^2 (G, B), \E)$ onto the set
  of complete, relatively continuous subspaces of~$\E$.  Its inverse is the
  map $\RC \mapsto \F (\E, \RC)$.

  A concrete Hilbert module~$\F$ is essential if and only if\/ $\RC_\F$ is
  dense.

  Isomorphism classes of Hilbert modules over $\Cred (G, B)$ correspond
  bijectively to isomorphism classes of continuously square-integrable Hilbert
  \Mpn{B,G}modules.
\end{theorem}

\begin{proof}
  Let~$\F$ be a concrete Hilbert \Mpn{A}module.  It is evident that $\RC_\F$
  is complete and relatively continuous.  Proposition~\ref{pro:RC_dense}
  asserts $\F (\E, \RC_\F) = \F$.  Conversely, let $\RC \subseteq \E$ be
  complete and relatively continuous.  Define $\F \defeq \F (\E, \RC)$.  Then
  $\RC \subseteq \RC_\F$.  We claim that $\RC = \RC_\F$.

  Let $x \in \RC_\F$, we want to show that $x \in \RC$.  Let $(u_j)_{j \in J}$
  be an approximate identity as in Lemma~\ref{lem:apprid}.  Since~$\F$ is the
  closure of $\KET {\RC}$, there is a sequence $(x_n) \in \RC$ with $\lim {}
  \KET {x_n} = \KET {x}$ in operator norm.
  Equation~\eqref{eq:right_estimate_si_KET} implies that
  $$
  \lim_{n \to \infty} {} \norm {x_n \ast u_j - x \ast u_j}_\si = 0
  $$
  for all $j \in J$.  Hence $x \ast u_j \in \RC$ because~$\RC$ is complete.
  Since $(u_j)$ is an approximate identity for~$A$, we have $\xi \cdot u_j \to
  \xi$ for all elements~$\xi$ of a Hilbert \Mpn{A}module.  In particular,
  $\KET {x} \cdot u_j = \KET {x \ast u_j}$ converges towards $\KET {x}$.
  Together with~\eqref{eq:apprid_right_module}, this means that $x \ast u_j
  \to x$ in the norm $\norm {\blank}_\si$.  Hence $x \in \RC$.  This proves
  that $\RC_\F = \RC$.

  If~$\F$ is essential, then $\RC_\F^0$ is dense in~$\E$.  Hence $\RC_\F$ is
  dense in~$\E$.  Conversely, if $\RC_\F$ is dense, then~$\F$ is essential by
  Proposition~\ref{pro:rel_cont}.  The last assertion of the theorem follows
  from Theorem~\ref{the:concrete_Hilbert}.
\end{proof}

\begin{corollary}  \label{cor:complete}
  Let $\RC \subseteq \E$ be relatively continuous.  Then the completion
  of\/~$\RC$ equals $\RC_{\F (\E, \RC)}$.  Thus the completion of\/~$\RC$ is
  still relatively continuous.
\end{corollary}

\begin{proof}
  The bijection between complete, relatively continuous subspaces and concrete
  Hilbert \Mpn{A}modules in Theorem~\ref{the:RC_characterize} preserves
  inclusions.  Hence $\RC_{\F (\E, \RC)}$ is the smallest complete subspace
  containing~$\RC$.
\end{proof}

\begin{corollary}  \label{cor:RC_G_B_module}
  Let $\RC \subseteq \E$ be a complete, relatively continuous subspace.
  Equip~$\RC$ with the norm $\norm {\blank}_\si$.  The subspace $\RC \subseteq
  \E$ is \Mpn{G}invariant, the action of~$G$ on~$\RC$ is continuous.
  Furthermore, $\RC$ is an essential right \Mpn{B}module, that is, $\RC \cdot
  B = \RC$.
\end{corollary}

\begin{proof}
  By Theorem~\ref{the:RC_characterize}, we have $\RC = \RC_\F$ for a concrete
  Hilbert module~$\F$ over~$A$.  Let $x \in \RC$.  By Cohen's Factorization
  Theorem, the map $a \mapsto \KET {x} \cdot a$ extends to a linear operator
  $\Mult (A) \to \F$ that is continuous with respect to the strict topology on
  $\Mult (A)$ and the norm topology on~$\F$.  We have $x \ast a \in \RC_\F$
  whenever there is $x \ast a \in \E_\si$ with $\KET {x \ast a} = \KET {x}
  \cdot a$ and $a \in \Mult (A)$.  Using~\eqref{eq:KET_B}
  and~\eqref{eq:KET_G}, we conclude that $x \cdot b, \gamma_g (b) \in \RC$ for
  all $b \in B$, $g \in G$.  Furthermore, we have norm estimates
  \eqref{eq:KET_B_est} and~\eqref{eq:KET_G_est}.  Since the map $g \mapsto
  \rho_g$ is strictly continuous, we have $\KET {x} \rho_g \to \KET {x}$ for
  $g \to 1$.  Therefore, the action of~$G$ on~$\RC$ is continuous.  Similarly,
  if $(u_i)$ is an approximate identity of~$B$, then $x \cdot u_i \to x$ in
  the norm $\norm {\blank}_\si$.  Hence Cohen's Factorization Theorem yields
  $\RC \cdot B = \RC$.
\end{proof}

Since~\eqref{eq:BRAKET} describes $\BRAKET {x} {y}$ explicitly, we may be able
to prove that $\BRAKET {x} {y} \in \Cred (G, B)$ without showing $x, y \in
\E_\si$.  For instance, it may happen that $\BRAKET {x} {y} \in C_c (G, B)$
for all $x, y \in \RC$.  If~$\RC$ is dense in~$\E$, then this implies $\RC
\subseteq \E_\si$:

\begin{proposition} \label{pro:si_automatic}
  Let $\RC \subseteq \E$ be a \emph{dense} subspace with $\BRAKET {x} {y} \in
  \Cred (G, B)$ for all $x, y \in \RC$.  Then $\RC \subseteq \E_\si$, so
  that~$\RC$ is relatively continuous.
\end{proposition}

\begin{proof}
  Fix $x \in \RC$.  If $f \in C_c (G, B)$, then $\KET {x} f$ is well-defined
  and
  $$
  \norm {\KET {x} f} =
  \norm {\braket {\KET {x} f} {\KET {x} f} }^{1/2} =
  \norm {\braket {f} {\BRAKET {x} {x} f}}^{1/2} \le
  \norm {f} \cdot \norm {\BRAKET {x} {x}}^{1/2}.
  $$
  Hence $\KET {x} \colon C_c (G, B) \to \E$ extends to a bounded operator
  $\KET {x} \colon L^2 (G, B) \to \E$.  The problem is to show that $\KET {x}$
  is adjointable.  This means that the set
  $$
  \E_0 \defeq \{ \xi \in \E \mid
  \text{$\exists f \in L^2 (G, B)$ with $\braket {\xi} {\KET {x} f_2} =
  \braket {f} {f_2}$ for all $f_2 \in C_c (G, B)$} \}
  $$
  equals~$\E$.  Since~$\E_0$ is a closed subspace, it suffices to prove
  that~$\E_0$ is dense.  If $y \in \RC$, $f_1, f_2 \in C_c (G, B)$, then
  $\braket {\KET {y} f_1} {\KET {x} f_2} = \braket {\BRAKET {x} {y} f_1}
  {f_2}$ and hence $\KET {y} f_1 \in \E_0$.  Elements of this form exhaust
  $\RC \ast C_c (G, B)$ by~\eqref{eq:right_module}.  Since~$\RC$ is dense
  in~$\E$, \eqref{eq:apprid_right_module} yields that $\RC \ast C_c (G, B)
  \supseteq \E_0$ is dense in~$\E$ as desired.
\end{proof}

\begin{definition}  \label{def:continuously_si}
  Let $(\E, \RC)$ and $(\E', \RC')$ be continuously square-integrable Hilbert
  modules.  We call $T \in \Adj (\E, \E')$ \emph{\Mpn{\RC}continuous} iff $T
  (\RC) \subseteq \RC'$ and $T^\ast (\RC') \subseteq \RC$.

  The \emph{generalized fixed point algebra} $\Fix (\E, \RC)$ is the closed
  linear span of $\KET {\RC} \BRA {\RC}$ in $\AdjG (\E)$.
\end{definition}

\begin{theorem}  \label{the:category_equivalence}
  Let $(\E, \RC)$ be a continuously square-integrable Hilbert \Mp{B,G}module
  and let $\F \defeq \F (\E, \RC)$.  There is a canonical, injective, strictly
  continuous \Mpn{\ast}homomorphism $\phi \colon \Adj (\F) \to \AdjG (\E)$,
  whose range is the space of \Mpn{\RC}continuous operators.  It maps $\Comp
  (\F)$ isometrically onto $\Fix (\E, \RC)$.

  The \Cstar{}algebra $\Fix (\E, \RC)$ is Morita-Rieffel equivalent to an
  ideal in $\Cred (G, B)$, namely the closed linear span of $\BRAKET {\RC}
  {\RC} \subseteq \Cred (G, B)$.

  The generalized fixed point algebra is the closed linear span of the
  operators $\int_G \gamma_g (x) \,dg$ with $x = \ket {\xi} \bra {\eta}$,
  $\xi, \eta \in \RC$.
\end{theorem}

\begin{proof}
  Since $\KET {\RC}$ is dense in~$\F$, we conclude that $\BRAKET {\RC} {\RC}$
  is dense in $\F^\ast \F$ and that $\KET {\RC} \BRA {\RC}$ is dense in $\F
  \F^\ast$.  Moreover, \eqref{eq:KET_compose} yields that the space~$M$
  defined in Theorem~\ref{the:concrete_morphism} equals the space of
  \Mpn{\RC}continuous operators.  Hence the assertions of the first paragraph
  follow from Theorem~\ref{the:concrete_morphism} if we take the
  homomorphism~$\phi$ defined there.  Since $\F A \subseteq \F$, the closed
  linear span~$J$ of $\F^\ast \F$ is an ideal in~$A$.  We may view~$\F$ as an
  imprimitivity bimodule for~$J$ and $\Comp (\F) \cong \Fix (\E, \RC)$.  That
  is, $\Fix (\E, \RC)$ and~$J$ are Morita-Rieffel equivalent.  The last
  assertion follows immediately from~\eqref{eq:KETBRA} and the definition of
  $\Fix (\E, \RC)$.
\end{proof}

Theorem~\ref{the:category_equivalence} implies that $(\E, \RC) \mapsto \F (\E,
\RC)$ is an equivalence between the \Cstar{}categories of continuously
square-integrable Hilbert \Mp{B,G}modules and Hilbert modules over $\Cred (G,
B)$, if we take \Mpn{\RC}continuous adjointable operators and adjointable
operators as morphisms, respectively.

Let~$\E$ be a square-integrable Hilbert \Mp{B,G}module.  It is an important
question whether there is a canonical choice for a dense, complete, relatively
continuous subset $\RC \subseteq \E$.  If~$B$ is proper, then there is one and
only one such~$\RC$.  This is the strongest sense in which~$\RC$ may be
canonical.  More generally, canonical should mean that we can single out a
specific subspace~$\RC$ using only that~$\E$ is a Hilbert \Mp{B,G}module.
Hence if $u \colon \E \to \E$ is an equivariant unitary, then $u (\RC) = \RC$
because~$u$ preserves the Hilbert \Mp{B,G}module structure.  This is
equivalent to $u \circ \F (\E, \RC) = \F (\E, \RC)$.  By
Corollary~\ref{cor:ideal} this happens iff all operators in $\AdjG (\E)$ are
\Mpn{\RC}continuous iff $\Fix (\E, \RC)$ is an ideal in $\AdjG (\E)$.

Unfortunately, $\AdjG (\E)$ frequently is so big that no ideal of it qualifies
as a generalized fixed point algebra.  For instance, if $B = \C$, then $\AdjG
(\E)$ will be a commutant of a group action on a Hilbert space and thus a von
Neumann algebra.  Hence there can be no canonical choice for~$\RC$ in this
case.  However, this does not yet create a very serious lack of uniqueness.
If $u \in \AdjG (\E)$ is unitary, then $(\E, \RC)$ and $(\E, u(\RC))$
correspond to two representations of the same abstract Hilbert module over
$\Cred (G, B)$ and hence give rise to isomorphic generalized fixed point
algebras.

\section{Constructions with relatively continuous subsets}
\label{sec:constructions}

As a preparation for Theorem~\ref{the:RC_tensor} and as an important special
case, we discuss Hilbert modules over $\Cred (G, B)$ of the form $\Cred (G,
\E)$.  Kasparov~\cite{kasparov:novikov} defines $\Cred (G, \E)$ as a
completion of $C_c (G, \E)$ with respect to a certain pre-Hilbert module
structure over $C_c (G, B)$.  An equivalent definition is
\begin{equation}  \label{eq:def_CredGE}
  \Cred (G, \E) \defeq \E \otimes_B \Cred (G, B),
\end{equation}
where we use the canonical map $B \to \Mult \bigl( \Cred (G, B) \bigr)$ to
form the tensor product.

In our framework, $\Cred (G, \E)$ arises as follows.  The subspace $C_c (G,
\E) \subseteq L^2 (G, \E)$ is dense and relatively continuous.
Equation~\eqref{eq:BRAKET} yields $\BRAKET {\xi} {\eta} \in C_c (G, B)$ for
all $\xi, \eta \in C_c (G, \E)$.  Hence $C_c (G, \E) \subseteq L^2 (G,
\E)_\si$ by Proposition~\ref{pro:si_automatic}.  We claim that
$$
\F \bigl( L^2 (G, \E), C_c (G, \E) \bigr) \cong \Cred (G, \E).
$$
To verify this, we generalize~\eqref{eq:rho_B} and define
\begin{equation}  \label{eq:rho_E}
  \rho_\eta \in \AdjG \bigl( L^2 (G, B), L^2 (G, \E) \bigr),
  \qquad
  (\rho_\eta f) (g) \defeq \gamma_g (\eta) \cdot f(g),
\end{equation}
for all $\eta \in \E$.  Equations \eqref{eq:KET_check}
and~\eqref{eq:KET_compose} yield
$$
\rho_\eta \circ \rho_K =
\rho_\eta \circ \KET {\check{K}} =
\KET {\rho_\eta (\check{K})}
$$
for all $\eta \in \E$, $K \in C_c (G, B)$.  It follows that $\rho (\E) \circ
\rho \bigl( C_c (G, B) \bigr)$ and $\KET {C_c (G, \E)}$ have the same closed
linear span in $\AdjG \bigl( L^2 (G, B), L^2 (G, \E) \bigr)$.
By~\eqref{eq:def_tensor}, the map
$$
\E \otimes_B \Cred (G, B) \to \F \bigl( L^2 (G, \E), C_c (G, \E) \bigr),
\qquad 
\eta \otimes K \mapsto \rho_\eta \circ \rho_K,
$$
is an isomorphism of Hilbert modules over $\Cred (G, B)$.

Consider the following situation.  Let $A$ and~$B$ be \Mpn{G}\Cstar{}algebras,
let $\E_1$ and $\E_2$ be \Mpn{G}equivariant Hilbert modules over $A$ and~$B$,
respectively, and let $\phi \colon A \to \Adj (\E_2)$ be an equivariant,
essential \Mpn{\ast}homomorphism.  The map~$\phi$ induces an essential
\Mpn{\ast}homomorphism $\Cred (G, A) \to \Adj \bigl( \Cred (G, \E_2) \bigr)$.

\begin{theorem}  \label{the:RC_tensor}
  Let $\RC_1 \subseteq \E_1$ be a (dense) relatively continuous subspace.  Let
  $\RC_{12}$ be the image of $\RC_1 \otimes^\alg \E_2 \subseteq \E_1
  \otimes^\alg \E_2$ under the canonical map to $\E_{12} \defeq \E_1 \otimes_A
  \E_2$.  Then $\RC_{12} \subseteq \E_{12}$ is (dense and) relatively
  continuous.  We have
  \begin{equation}  \label{eq:RC_tensor}
    \F (\E_{12}, \RC_{12}) \cong
    \F (\E_1, \RC_1) \otimes_{\Cred (G, A)} \Cred (G, \E_2).
  \end{equation}
\end{theorem}

\begin{proof}
  Let $\xi \in \RC_1$, $\eta \in \E_2$.  Since~$\phi$ is essential, $L^2 (G,
  A) \otimes_A \E_2 \cong L^2 (G, \E_2)$.  Hence $\KET {\xi} \otimes_A
  \ID_{\E_2} \in \AdjG (L^2 (G, \E_2), \E_{12})$.  The same simple computation
  that yields~\eqref{eq:KET_B} shows that
  $$
  \KET {\xi \otimes \eta} =
  (\KET {\xi} \otimes \ID_{\E_2}) \circ \rho_\eta,
  $$
  where $\rho_\eta$ is defined by~\eqref{eq:rho_E}.  As a result, $\RC_{12}
  \subseteq (\E_{12})_\si$ and
  $$
  \KET {\RC_{12}} =
  \KET {\RC_1 \otimes \E_2} =
  (\KET {\RC_1} \otimes \ID_{\E_2}) \circ \rho (\E_2).
  $$
  By definition, $\F_{12} \defeq \F (\E_{12}, \RC_{12})$ is the closed linear
  span of $\KET {\RC_{12}} \circ \Cred (G, B)$.  The discussion of $\Cred (G,
  \E)$ above shows that the closed linear span of $\rho (\E_2) \circ \Cred (G,
  B)$ equals $\Cred (G, \E_2)$.  Hence $\F_{12}$ equals the closed linear span
  of $(\F_1 \otimes \ID_{\E_2}) \circ \Cred (G, \E_2)$.
  Equation~\eqref{eq:def_tensor} yields that the map
  $$
  \F_1 \otimes_{\Cred (G, A)} \Cred (G, \E_2) \to \F_{12},
  \qquad
  \xi \otimes \eta \mapsto (\xi \otimes_A \ID_{\E_2}) \circ \eta
  $$
  is an isometry of Hilbert modules over $\Cred (G, B)$.  Hence $\RC_{12}$ is
  relatively continuous and satisfies~\eqref{eq:RC_tensor}.
\end{proof}

A consequence of the proof (or of~\eqref{eq:RC_tensor} and
Corollary~\ref{cor:complete}) is that the construction $\RC_1 \mapsto
\RC_{12}$ is compatible with completions.  That is, if the completions of
$\RC_1$ and $\RC_1'$ are equal, then the same holds for $\RC_{12}$ and
$\RC_{12}'$.

We consider some important special cases of Theorem~\ref{the:RC_tensor}.

\begin{corollary}  \label{cor:RC_transport}
  Let $(A, \RC)$ be a continuously square-integrable \Cstar{}algebra and let
  $\phi \colon A \to \Adj (\E)$ be an equivariant, essential
  \Mpn{\ast}homomorphism.

  Then $\RC (\E) \subseteq \E$ is a dense, relatively continuous subset, and
  $$
  \F \bigl(\E, \RC (\E) \bigr) \cong
  \F (A, \RC) \otimes_{\Cred (G,A)} \Cred (G, \E).
  $$
\end{corollary}

\begin{proof}
  The isomorphism $A \otimes_A \E \cong \E$ maps $\RC \otimes \E$ onto the
  linear span of $\RC (\E)$.
\end{proof}

In particular, if $\Comp (\E)$ is continuously square-integrable, so is~$\E$.
Conversely:

\begin{corollary}  \label{cor:RC_comp}
  Let $(\E, \RC)$ be a continuously square-integrable Hilbert \Mp{B,G}module.
  Then the linear span of
  $$
  \ket {\RC} \bra {\E} \defeq
  \{ \ket{\xi} \bra{\eta} \mid \xi \in \RC,\ \eta \in \E \}
  $$
  is a dense, relatively continuous subspace of\/ $\Comp (\E)$.  We have
  $$
  \F \bigl( \Comp (\E), \ket {\RC} \bra {\E} \bigr) \cong
  \F (\E, \RC) \otimes_{\Cred (G, B)} \Cred (G, \E^\ast).
  $$
\end{corollary}

\begin{proof}
  The assertion follows from \eqref{eq:Comp_ketbra} and
  Theorem~\ref{the:RC_tensor}.
\end{proof}

Let $S(\E)$ and $S \Comp (\E)$ be the sets of all dense, complete, relatively
continuous subspaces of $\E$ and $\Comp (\E)$, respectively.  Corollaries
\ref{cor:RC_transport} and~\ref{cor:RC_comp} give rise to maps
$$
i \colon S(\E) \to S \Comp (\E),
\qquad
j \colon S \Comp (\E) \to S(\E).
$$
We analyze whether these two maps are inverse to each other.  Recall that a
group~$G$ is \emph{exact} if and only if
$$
\Cred (G, I) = \ker \bigl( \Cred (G, B) \to \Cred (G, B/I) \bigr)
$$
whenever $I \subseteq B$ is an invariant closed ideal in a
\Mpn{G}\Cstar{}algebra~$B$.

\begin{theorem}  \label{the:EKE}
  For all $G$ and~$B$ and all Hilbert \Mp{B,G}modules~$\E$, the composition $i
  \circ j$ is the identity map on $S \Comp (\E)$ and $j \circ i
  (\RC) \subseteq \RC$ for all $\RC \in S (\E)$.

  If the group~$G$ is exact, then $j \circ i$ is the identity map on $S(\E)$.
  Conversely, if~$G$ is not exact, then there are Hilbert modules~$\E$ for
  which $j \circ i$ is not the identity map.
\end{theorem}

\begin{proof}
  Let $\RC \in S \Comp (\E)$.  Then $j (\RC)$ is the completion of $\RC (\E)
  \subseteq \E$.  Since~$i$ is compatible with completions, $i \circ j (\RC)$
  is the completion of $\ket{\RC (\E)} \bra {\E} = \RC \circ \ket {\E} \bra
  {\E}$.  The linear span of $\RC \circ \ket {\E} \bra {\E}$ is a dense
  subspace of~$\RC$ by Corollary~\ref{cor:RC_G_B_module}.  Hence the
  completion of $\RC \circ \ket {\E} \bra {\E}$ equals~$\RC$.  This proves $i
  \circ j = \ID$.

  Conversely, let $\RC \in S (\E)$ and $\F \defeq \F (\E, \RC)$.  Then $j
  \circ i (\RC)$ is the completion of $\ket {\RC} \bra {\E} (\E) = \RC \cdot
  \braket {\E} {\E}$.  Hence $\F \bigl(\E, ji (\RC) \bigr)$ is the closed
  linear span of $\F \cdot \braket {\E} {\E}$.  Let $I \subseteq B$ be the
  closed ideal generated by $\braket {\E} {\E}$.  Cohen's Factorization
  Theorem and Corollary~\ref{cor:RC_G_B_module} show that
  $$
  j \circ i (\RC) = \RC \cdot I \subseteq \RC,
  \qquad
  \F \bigl(\E, ji (\RC) \bigr) = \F \cdot I \subseteq \F.
  $$
  Let $J \subseteq \Cred (G, B)$ be the closed ideal generated by $\braket
  {\F} {\F}$.  If $J \subseteq \Cred (G, I)$, then we may view~$\F$ as a
  Hilbert module over $\Cred (G, I)$.  Hence $\F \cdot I = \F$.  Conversely,
  if $\F \cdot I = \F$, then $J \cdot I = J$ and hence $J \subseteq \Cred (G,
  B) \cdot I = \Cred (G, I)$.  Therefore,
  $$
  j \circ i (\RC) = \RC \iff J \subseteq \Cred (G, I).
  $$

  If $\xi, \eta \in \RC$, then $\BRAKET {\xi} {\eta} (g) \in I$ for all $g \in
  G$ by~\eqref{eq:BRAKET}.  Therefore, $\BRAKET {\xi} {\eta}$ is annihilated
  by the canonical map $\Cred (G, B) \to \Cred (G, B/I)$.  If~$G$ is exact,
  this implies that $\BRAKET {\xi} {\eta} \in \Cred (G, I)$.  Hence $J
  \subseteq \Cred (G, I)$ and thus $j \circ i (\RC) = \RC$.

  Suppose that~$G$ is not exact and that $I \subseteq B$ is an invariant ideal
  for which $\Cred (G, I)$ is strictly smaller than the kernel~$K$ of the map
  $\Cred (G, B) \to \Cred (G, B/I)$.  View~$K$ as a Hilbert module over $\Cred
  (G, B)$ and let $(\E, \RC)$ be the associated continuously square-integrable
  Hilbert module according to Theorem~\ref{the:RC_characterize}.  Then
  $\BRAKET {\RC} {\RC} \subseteq K$.  This implies $\braket {\xi} {\eta} \in
  I$ for all $\xi, \eta \in \RC$ by~\eqref{eq:BRAKET}.  Hence $\braket {\E}
  {\E} \subseteq I$.  However, $J = K$ is not contained in $\Cred (G, I)$.
  Hence $\RC \neq \RC \cdot I$.
\end{proof}

We remark that the identity $i \circ j = \ID$ is equivalent to the isomorphism
$$
\Cred (G, \E) \otimes_{\Cred (G, B)} \Cred (G, \E^\ast) \cong
\Cred (G, \E \otimes_B \E^\ast) \cong \Cred \bigl(G, \Comp (\E) \bigr).
$$

\section{Some counterexamples}
\label{sec:counter}

In this section, we consider a simple special case in which a complete
description of the square-integrable and continuously square-integrable
Hilbert modules is possible.  We assume that $B = \C$ and that~$G$ is Abelian,
\Mpn{\sigma}compact, and metrizable, but not compact.  Hence the Pontrjagin
dual~$\hat{G}$ of~$G$ is not discrete.  For instance, we may take $G = \Z^n$
for some $n \in \N \setminus \{0\}$.

Since $\Cred (G, B) \cong C_0 (\hat{G})$, countably generated Hilbert modules
over $\Cred (G, B)$ correspond to continuous fields of separable Hilbert
spaces over~$\hat{G}$.  Let $(\Hils_x)_{x \in \hat{G}}$ be a continuous field
of Hilbert spaces.  The associated Hilbert module over $C_0 (\hat{G})$ is $C_0
\bigl(\hat{G}, (\Hils_x) \bigr)$, the space of continuous sections of
$(\Hils_x)$ vanishing at infinity.  The \Cstar{}algebra of compact operators
on this Hilbert module is isomorphic to $C_0 \bigl(\hat{G}, \Comp (\Hils_x)
\bigr)$, where $\bigl( \Comp (\Hils_x) \bigr)_{x \in \hat{G}}$ carries the
canonical bundle structure.

A countably generated Hilbert \Mpn{B,G}module is nothing but a representation
of~$G$ on a separable Hilbert space.  By the Equivariant Stabilization
Theorem, a \Mpn{G}Hilbert space is square-integrable if and only if it is a
direct summand in $(L^2 G)^\infty \cong L^2 (\hat{G}, dx)^\infty$, where~$dx$
denotes the Haar measure on~$\hat{G}$.  Therefore, a \Mpn{G}Hilbert space is
square-integrable iff it is equivalent to a Hilbert space of square-integrable
sections of some measurable field of Hilbert spaces over~$\hat{G}$, equipped
with the canonical representation of~$G$ by pointwise multiplication.  Two
measurable fields yield equivalent representations of~$G$ iff they are
isomorphic outside a set of Haar measure zero.

Measurable fields of Hilbert spaces are classified by the \emph{dimension
function}
$$
d \colon \hat{G} \to \exN \defeq \N \cup \{\infty, 0\}
$$
that associates to $x\in \hat{G}$ the dimension of the fiber over~$x$.  The
function~$d$ is measurable, and any measurable function arises as the
dimension function of a measurable field.  We say that an assertion holds
a.e.\ (almost everywhere) iff it holds outside a set of measure zero.  Two
measurable fields are isomorphic a.e.\ if and only if the dimension functions
agree a.e..  Hence isomorphism classes of square-integrable, separable
\Mpn{G}Hilbert spaces correspond to a.e.-equality classes of measurable
functions $\hat{G} \to \exN$.

Let $(\Hils_x)_{x \in \hat{G}}$ be a continuous field of Hilbert spaces
over~$\hat{G}$.  If we view $(\Hils_x)$ as a Hilbert module over $\Cred G$ and
apply the functor $\blank \otimes_{\Cred G} L^2G$, we get the Hilbert space of
square-integrable sections of $(\Hils_x)_{x \in \hat{G}}$ with the
representation of~$G$ by pointwise multiplication.  Hence the functor $\blank
\otimes_{\Cred G} L^2G$ forgets everything about the field $(\Hils_x)$ except
the a.e.-equality class of its dimension function.

The dimension function of a continuous field of Hilbert spaces is
automatically lower semi-continuous.  Therefore, if $d \colon \hat{G} \to
\exN$ is not equal a.e.\ to a lower semi-continuous function, then the
corresponding square-integrable representation cannot come from a Hilbert
module over $\Cred G$.  To construct examples of such measurable functions,
let~$d$ be the characteristic function of a compact subset $K \subseteq
\hat{G}$.  Suppose that $d' \colon \hat{G} \to \exN$ is lower semi-continuous
and that $d' \le d$ a.e..  It follows that $d' \le 1$ and that $d' = 0$ on the
open set $\hat {G} \setminus K$.  Thus~$d'$ is the characteristic function of
an open subset $U \subseteq K$.  If~$K$ is a compact set with non-zero Haar
measure and empty interior, then $d' = 0$ is the only lower semi-continuous
dimension function with $d' \le d$ a.e..  Nevertheless, $d' \neq d$ a.e..

The square-integrable \Mpn{G}Hilbert space associated to~$d$ is $L^2 (K, dx)$,
on which~$G$ acts by pointwise multiplication.  Suppose that $\RC \subseteq
L^2 (K, dx)$ is relatively continuous and complete.  Let $\Hils \subseteq L^2
(K)$ be the closure of~$\RC$.  Then $(\Hils, \RC)$ is continuously
square-integrable.  Therefore, the dimension function of~$\Hils$ is lower
semi-continuous.  By construction of~$d$ this implies that $\Hils = \{0\}$.
Consequently, $\{0\}$ is the only relatively continuous subset of the
square-integrable \Mpn{G}Hilbert space $L^2 (K, dx)$.

Conversely, if the dimension function of a separable \Mpn{G}Hilbert space is
lower semi-continuous, there is a dense, relatively continuous subspace.  The
proof is left to the reader.  However, this subspace is never unique.  Even
more, there are many Hilbert modules~$\F$ over $\Cred G$ for which $\F
\otimes_{\Cred G} L^2 G \cong L^2 G$.  The most obvious source of
non-uniqueness is modification on a set of measure zero.  Let $S \subseteq
\hat{G}$ be a closed subset of measure zero (for instance, a finite subset).
The ideal
$$
I_S = \{ f \in C_0 (\hat{G}) \mid f|_S = 0 \} \subseteq C_0 (\hat{G})
$$
may be viewed as a Hilbert module over $C_0 (\hat{G})$ and thus as a
continuous field of Hilbert spaces over~$\hat{G}$.  Its dimension function is
the characteristic function of $\hat{G} \setminus S$ and hence equal to~$1$
a.e..  Thus
$I_S \otimes_{\Cred G} L^2 G \cong L^2 G$.
The generalized fixed point algebra in this example is $I_S$.  Hence a
generalized fixed point algebra for $L^2 G$ need not be isomorphic to $C_0
(\hat{G})$, not even Morita-Rieffel equivalent to $C_0 (\hat{G})$.

The lack of uniqueness observed above can be overcome by restricting attention
to \emph{maximal} relatively continuous subsets, that is, relatively
continuous subsets that are not contained in any larger relatively continuous
subset.  Since $I_S \subseteq C_0 (\hat{G})$, the subspace $\RC_{I_S}$ cannot
be maximal.  However, even if we insist on maximality, we do not obtain
uniqueness of the generalized fixed point algebra, because a continuous field
is not yet determined by its dimension function.

If $(\Hils_x)_{x \in \hat{G}}$ is a continuous field of Hilbert spaces with
$\dim \Hils_x = n$ for all $x \in \hat{G}$, then $(\Hils_x)$ ``is'' an
\Mpn{n}dimensional complex vector bundle over~$\hat{G}$.  That is, there is an
\Mpn{n}dimensional complex vector bundle~$E$ over~$\hat{G}$ such that the
space of continuous sections of $(\Hils_x)$ is isomorphic to $C_0 (\hat{G},
E)$ as a module over $C_0 (\hat{G})$.  The corresponding square-integrable
representation of~$G$ is $L^2 (G)^n$ because the dimension function is
constant.

We denote the \Mpn{n}dimensional trivial vector bundle by~$\C^n$.
The generalized fixed point algebra associated to a vector bundle $E \to
\hat{G}$ is $C_0 \bigl(\hat{G}, \End (E) \bigr)$.  Especially,
$$
C_0 \bigl(\hat{G}, \End (\C^n) \bigr) \cong C_0 (\hat{G}, \Mat_n).
$$
If~$E$ is a line bundle, then $\End (E) \cong E \otimes E^\ast$ is trivial.
Hence the generalized fixed point algebras associated to vector bundles $E$
and~$E'$ are isomorphic if $E' \cong E \otimes L$ for a complex line
bundle~$L$.  The converse also holds: If the generalized fixed point algebras
are isomorphic, then $E$ and~$E'$ differ by tensoring with a line bundle.  We
leave the proof as an exercise in vector bundle theory for the interested
reader.

In particular, the generalized fixed point algebra equals $C_0 (\hat{G},
\Mat_n)$ if and only if $E \cong \C^n \otimes L = L \oplus L \oplus \dots
\oplus L$ is a direct sum of~$n$ copies of the same line bundle.  This can
be shown easily by observing how the matrix units in~$\Mat_n$ operate on~$E$.
Hence if $E \oplus \C$ is of this form, then~$E$ has to be trivial.  As a
result, whenever there is a non-trivial vector bundle over~$\hat{G}$, we can
find one for which the generalized fixed point algebra is not isomorphic to
$C_0 (\hat{G}, \Mat_n)$.  However, the resulting generalized fixed point
algebras are always Morita-Rieffel equivalent to $C_0 (\hat{G})$.

We claim that the relatively continuous subset~$\RC$ of $L^2 (G)^n$ associated
to a vector bundle $(\Hils_x)$ over~$\hat{G}$ is always maximal.  Hence we
cannot rule out the lack of uniqueness of $\Fix (L^2 (G)^n, \RC)$ by requiring
maximality.  If $\RC \subseteq \RC'$ and~$\RC'$ is relatively continuous, then
$\F (L^2 (G)^n, \RC')$ corresponds to a continuous field of Hilbert spaces
$(\Hils_x')$ over~$\hat{G}$.  We have $\Hils_x \subseteq \Hils'_x$ for all~$x$
and $\dim \Hils_x' \le n$ outside a set of measure zero.  Lower
semi-continuity implies that $\dim \Hils'_x = \dim \Hils_x = n$ for all~$x$.
Therefore, $(\Hils_x') = (\Hils_x)$ and hence $\RC' = \RC$.

Finally, we claim that no non-zero continuously square-integrable
\Mpn{G}Hilbert space is ideal.  Let $(\Hils_s)_{s \in \hat{G}}$ be a non-zero
continuous field of Hilbert spaces over~$\hat{G}$.  Choose a non-zero
continuous section $f \in C_0 \bigl( \hat{G}, (\Hils_s) \bigr)$.  Let $U
\subseteq \hat{G}$ be an open set with $f (s) \neq 0$ for all $s \in U$.
There is a bounded, positive, measurable function $\phi \colon \hat{G} \to \C$
whose restriction to~$U$ is not equal a.e.\ to a continuous function.  Hence
$\phi^2 \cdot \norm{f}^2$ is not equal a.e.\ to a continuous function, so that
the square-integrable section $\phi \cdot f$ of $(\Hils_x)$ is not continuous.
Therefore, the operator of pointwise multiplication by~$\phi$ is not
\Mpn{\RC}continuous, although it is equivariant and adjointable.

\section{Proper coefficients}
\label{sec:proper}

The last section shows that there are significant differences between
continuously square-integrable, square-integrable, and arbitrary equivariant
Hilbert modules for $B = \C$.  Nevertheless, if the group action on~$B$ is
``sufficiently proper'', these differences disappear.  That is, any Hilbert
\Mp{B,G}module is square-integrable and contains a unique dense, complete,
relatively continuous subspace.  This happens if~$B$ is proper in Kasparov's
sense and, more generally, if the induced group action on the (not necessarily
separated) spectrum of~$B$ is proper.

\begin{definition}  \label{def:proper}
  Let~$X$ be a not necessarily separated topological space.  Let~$G$ be a
  locally compact group and let $G \times X \to X$ be a continuous action
  of~$G$ on~$X$.  We call~$X$ a \emph{proper \Mpn{G}space} iff for all $x,y
  \in X$ there are neighborhoods $U_x$ and~$U_y$ of $x$ and~$y$ in~$X$ such
  that the set
  $$
  \{ g \in G \mid g(U_x) \cap U_y \neq \emptyset \} \subseteq G
  $$
  is relatively compact.
\end{definition}

We call $K \subseteq X$ \emph{quasi-compact} iff any open covering of~$K$ has
a finite subcovering and \emph{relatively quasi-compact} iff~$K$ is contained
in a quasi-compact subset of~$X$.

The following lemma shows that Definition~\ref{def:proper} contains the usual
definition of proper actions on separated, locally compact spaces.

\begin{lemma}  \label{lem:proper}
  Let~$X$ be a not necessarily separated, proper \Mpn{G}space.  Let $K, L
  \subseteq X$ be relatively quasi-compact.  Then there are open neighborhoods
  $U_K$ and~$U_L$ of $K$ and~$L$, respectively, such that the set
  $$
  \{ g \in G \mid g(U_K) \cap U_L \neq \emptyset \} \subseteq G
  $$
  is relatively compact.
\end{lemma}

\begin{proof}
  The proof is a straightforward exercise in topology.  We may assume without
  loss of generality that $K$ and~$L$ are quasi-compact, not just relatively
  quasi-compact.  For $x \in K$, $y \in L$, there are open neighborhoods
  $U^{xy}_x$, $U^{xy}_y$ of $x$ and~$y$ such that $g(U^{xy}_x) \cap U^{xy}_y =
  \emptyset$ for all~$g$ outside a compact subset of~$G$, because the action
  on~$X$ is proper.  By quasi-compactness, for fixed~$x$ finitely many of the
  open sets $U_y^{xy}$ cover~$L$.  Let $U_L^x$ be their union and let $U_x^x$
  be the intersection of the corresponding $U^{xy}_x$.  Then $(U_x^x)_{x \in
  K}$ is an open covering of~$K$.  Finitely many of these sets suffice to
  cover~$K$.  Let~$U_K$ be their union and let $U_L$ be the intersection of
  the corresponding open neighborhoods~$U_L^x$ of~$L$.  These sets have the
  desired properties.
\end{proof}

Let~$B$ be a \Mpn{G}\Cstar{}algebra and let~$P$ be its primitive ideal space,
equipped with the Jacobson topology and the continuous action of~$G$ defined
by $g \cdot \prid \defeq \{ \beta_g (b) \mid b \in \prid \}$ for $g \in G$ and
$\prid \in P$.  It makes no difference to use the space of irreducible
representations of~$B$ instead because we only use the lattice of open subsets
of~$P$.

\begin{definition}  \label{def:spec_proper}
  A \Mpn{G}\Cstar{}algebra is called \emph{spectrally proper} iff its
  primitive ideal space is a proper \Mpn{G}space.
\end{definition}

We claim that proper \Mpn{G}\Cstar{}algebras are spectrally proper.  Let~$X$
be a proper, locally compact \Mpn{G}space.  By the Dauns-Hoffmann theorem, the
center of $\Mult (B)$ is isomorphic to $C_b (P)$.  It follows that essential
\Mpn{\ast}homomorphisms from $C_0 (X)$ to the center of $\Mult (B)$ correspond
to continuous maps $P \to X$ \cite{Nilsen:bundles}.  As a result, $B$ is a
proper \Mpn{G}\Cstar{}algebra iff there is an equivariant, continuous map $P
\to X$ for a separated, locally compact, proper \Mpn{G}space~$X$.  This
implies that~$P$ is proper, that is, $B$ is spectrally proper.

We recall some well-known facts about the primitive ideal space to fix our
notation.  If $b \in B$, $\prid \in P$, let $b_\prid$ be the image of~$b$ in
the quotient $B/\prid$.  Open subsets $U \subseteq P$ correspond to closed
ideals in~$B$ via
$$
U \mapsto B_U \defeq \bigcap_{\prid \in P \setminus U} \prid =
\{ b \in B \mid
\text{$b_\prid = 0$ for all $\prid \in P \setminus U$} \}.
$$
We have $U_1 \subseteq U_2$ if and only if $B_{U_1} \subseteq B_{U_2}$.  If
$U_1, U_2 \subseteq P$ are relatively quasi-compact and open, the same holds
for $U_1 \cup U_2$.  Hence the family~$\mathcal{C}$ of all relatively
quasi-compact, open subsets of~$P$ is directed.  Therefore, the union
$$
B_\cs \defeq \bigcup_{U \in \mathcal{C}} B_U
$$
is a \Mpn{\ast}ideal in~$B$.  This ideal is dense in~$B$ because the sets
$$
U_{b,t} \defeq \{ \prid \in P \mid \norm {b_\prid} > t \}
$$
are open and relatively quasi-compact for all $b \in B$, $t > 0$.  Functional
calculus allows us to approximate~$b$ in norm by elements of $U_{b,t}$ with $t
> 0$.

\begin{theorem}  \label{the:proper_csi}
  Let~$B$ be a spectrally proper \Mpn{G}\Cstar{}algebra and let~$\E$ be a
  \Mpn{G}equivariant Hilbert module over~$B$.  Let $B_\cs$ be defined as above
  and let $\E_\cs \defeq \E \cdot B_\cs$.  Then $\BRAKET {\E_\cs} {\E_\cs}
  \subseteq C_c (G, B)$ and $\E_\cs$ is a dense, relatively continuous
  subspace of~$\E$.  In particular, $\E$ is square-integrable.

  The completion $\RC_0$ of $\E_\cs$ is the \emph{only} dense, complete,
  relatively continuous subspace of~$\E$.  We have $\xi \in \RC_0$ if and only
  if $\xi \in \E_\si$ and $\BRAKET {\xi} {\xi} \in \Cred (G, B)$.  Any
  relatively continuous subset of~$\E$ is contained in $\RC_0$.
\end{theorem}

\begin{proof}
  Let $U \subseteq P$ be open and relatively compact and let $\E_U \defeq \E
  \cdot B_U$.  By Cohen's Factorization Theorem, $\E_U$ is a closed linear
  subspace of~$\E$.  The subset
  $$
  V \defeq \{ g \in G \mid gU \cap U \neq \emptyset \} \subseteq G
  $$
  is open and relatively compact by Lemma~\ref{lem:proper}.  We may write
  elements of $\E_U$ in the form $\xi \cdot b$, $\eta \cdot c$ with $\xi, \eta
  \in \E$, $b,c \in B_U$.  Equation~\eqref{eq:BRAKET} yields
  $$
  \BRAKET {\xi \cdot b} {\eta \cdot c} (g) =
  b^\ast \cdot \braket {\xi} {\gamma_g (\eta)} \cdot \beta_g (c) \in
  B_U \cdot B \cdot \beta_g (B_U) \subseteq
  B_{g U \cap U}.
  $$
  Hence $\BRAKET {\E_U} {\E_U} \subseteq C_0 (V, B) \subseteq C_c (G, B)$.  It
  follows that $\BRAKET {\E_\cs} {\E_\cs} \subseteq C_c (G, B)$.  Since
  $B_\cs$ is dense in~$B$, $\E_\cs$ is dense in~$\E$.
  Proposition~\ref{pro:si_automatic} yields $\E_\cs \subseteq \E_\si$, so that
  $\E_\cs \subseteq \E$ is a dense, relatively continuous subspace and~$\E$ is
  square-integrable.

  Since $\norm {\BRAKET {\xi} {\xi} (g)} \le \norm {\xi}^2$, it follows that
  $\norm {\KET {\xi}} \le C_U \cdot \norm {\xi}$ for all $\xi \in \E_U$ with
  some $C_U > 0$.  Hence the norms $\norm {\blank}$ and $\norm {\blank}_\si$
  are equivalent on $\E_U$.

  Let $\RC \subseteq \E$ be a dense, complete, relatively continuous subspace.
  We claim that~$\RC$ contains $\E_\cs$, so that $\RC_0 \subseteq \RC$.
  Corollary~\ref{cor:RC_G_B_module} implies that $\RC \cdot B_U \subseteq
  \RC$.  Since~$\RC$ is dense in~$\E$, $\RC \cdot B_U$ is dense in $\E_U$ with
  respect to the norm $\norm {\blank}$ and hence also with respect to the norm
  $\norm {\blank}_\si$ because these two norms are equivalent on~$\E_U$.
  Since $\RC$ is complete, it follows that $\E_U \subseteq \RC$.  Hence
  $\E_\cs \subseteq \RC$.

  If $\xi \in \RC_0$, then $\xi \in \E_\si$ and $\BRAKET {\xi} {\xi} \in \Cred
  (G, B)$ by Corollary~\ref{cor:complete}.  Assume conversely that $\xi \in
  \E_\si$ and $\BRAKET {\xi} {\xi} \in \Cred (G, B)$.  We claim that $\xi \in
  \RC_0$.  Since $B_\cs \subset B$ is a dense \Mpn{\ast}ideal, there is an
  approximate identity $(u_i)_{i \in I}$ for~$B$ with $u_i \in B_\cs$ for all
  $i \in I$.  Thus $\xi \cdot u_i \in \E_\cs \subseteq \RC_0$ for all $i \in
  I$.  Let $\RC \subseteq \E_\si$ be the completion of $\{\xi\}$.
  Corollary~\ref{cor:RC_G_B_module} applied to~$\RC$ yields that $\xi \cdot
  u_i \to \xi$ in the norm $\norm {\blank}_\si$.  Hence $\xi \in \RC_0$ as
  asserted.

  Therefore, any relatively continuous subset $\RC \subseteq \E$ is contained
  in $\RC_0$.
\end{proof}

As a consequence, we obtain the following result of Kasparov and Skandalis.

\begin{corollary}  \label{cor:csi_category}
  Let~$B$ be a spectrally proper \Mpn{G}\Cstar{}algebra.  Then the functor
  $$
  \F \mapsto \F \otimes_{\Cred (G, B)} L^2 (G, B)
  $$
  is an equivalence between the \Cstar{}categories of Hilbert modules over
  $\Cred (G, B)$ and \Mpn{G}equivariant Hilbert modules over~$B$.  That is,
  any \Mpn{G}equivariant Hilbert module~$\E$ over~$B$ arises in this way for a
  unique Hilbert module~$\F$ over $\Cred (G, B)$, and the map $\Adj (\F) \to
  \AdjG (\E)$ is an isomorphism.
\end{corollary}

\begin{proof}
  Let~$\E$ be a Hilbert \Mp{B,G}module.  By Theorem~\ref{the:proper_csi},
  there is a unique dense, complete, relatively continuous subset $\RC
  \subseteq \E$.  Hence there is no difference between continuously
  square-integrable Hilbert \Mp{B,G}modules and Hilbert \Mp{B,G}modules.
  Theorem~\ref{the:RC_characterize} shows that isomorphism classes of Hilbert
  modules over $\Cred (G, B)$ and Hilbert \Mp{B,G}modules correspond to each
  other bijectively.  Since $\RC \subseteq \E$ is unique, we have $u (\RC) =
  \RC$ for all $u \in \AdjG (\E)$.  Hence $\AdjG (\E) \cong \Adj (\F)$ by
  Corollary~\ref{cor:ideal}.
\end{proof}

\bibliography{fixed}
\bibliographystyle{amsplain}

\end{document}